\documentclass{amsart}
\usepackage{amssymb,verbatim}
\usepackage{hyperref}

 \DeclareMathOperator{\diam}{diam}
\DeclareMathOperator{\E}{\mathbb{E}}
\DeclareMathOperator{\EE}{\mathbb{E}}
\DeclareMathOperator{\sgn}{sgn}

\newcommand{\M}{\mathcal{M}}
\newcommand{\N}{\mathcal{N}}

\newcommand{\F}{\mathcal{F}}
\newcommand{\R}{\mathbb{R}}
\newcommand{\B}{\mathcal B}
\newcommand{\Z}{\mathbb{Z}}
\newcommand{\G}{\mathcal{G}}

\newcommand{\Rad}{{\mathrm{\bf Rad}}}
\newcommand{\Lip}{\mathrm{Lip}}
\newcommand{\op}{\mathcal{E}}
\newcommand{\diag}{\mathrm{\bf diag}}
\newcommand{\dist}{\mathrm{dist}}
\newcommand{\sign}{\mathrm{sign}}
\newcommand{\Le}{\Biggl}
\newcommand{\Ri}{\Biggr}

\newcommand{\e}{\varepsilon}

\newtheorem{theorem}{Theorem}[section]
\newtheorem{lemma}[theorem]{Lemma}
\newtheorem{prop}[theorem]{Proposition}
\newtheorem{claim}[theorem]{Claim}

\newtheorem{corollary}[theorem]{Corollary}

\theoremstyle{definition}

\newtheorem{definition}[theorem]{{\it Definition}}

\newtheorem{remark}[theorem]{{\it Remark}}

\begin{document}
\title{Metric cotype}

\author{Manor Mendel}
\thanks{MM was partially supported by ISF grant no. 221/07,
BSF grant no. 2006009, and
a gift from Cisco research center}
\address {Computer Science Division,
The Open University of Israel,  
Israel}
\email{mendelma@gmail.com}. 
\author{Assaf Naor}
\address{Courant Institute of Mathematical Sciences, New York University, New York, NY}
\email{naor@cims.nyu.edu}

\begin{abstract}
We introduce the notion of {\em cotype of a metric space}, and
prove that for Banach spaces it coincides with the classical
notion of Rademacher cotype. This yields a concrete version of
Ribe's theorem, settling a long standing open problem in the
nonlinear theory of Banach spaces. We apply our results to
several problems in metric geometry. Namely, we use metric cotype
in the study of uniform and coarse embeddings, settling in
particular the problem of classifying when $L_p$ coarsely or
uniformly embeds into $L_q$. We also prove a nonlinear analog of
the Maurey-Pisier theorem, and use it to answer a question posed
by Arora, Lov\'asz, Newman, Rabani, Rabinovich and Vempala, and to
obtain quantitative bounds in a metric Ramsey theorem due to
Matou\v{s}ek.
 \end{abstract}

\maketitle

\section{Introduction}

In 1976 Ribe~\cite{Ribe76} (see
also~\cite{Ribe78}, \cite{HM82}, \cite{Bourgain87}, \cite{BL00}) proved that if $X$ and
$Y$ are uniformly homeomorphic Banach spaces then $X$ is finitely
representable in $Y$, and vice versa ($X$ is said to be finitely
representable in $Y$ if there exists a constant $K>0$ such that
any finite dimensional subspace of $X$ is $K$-isomorphic to a
subspace of $Y$). This theorem suggests that ``local properties"
of Banach spaces, i.e. properties whose definition involves
statements about finitely many vectors, have a purely metric
characterization. Finding explicit manifestations of this
phenomenon for specific local properties of Banach spaces (such as
type, cotype and super-reflexivity), has long been a major driving
force in the bi-Lipschitz theory of metric spaces (see Bourgain's
paper~\cite{Bourgain86-trees} for a discussion of this research
program). Indeed, as will become clear below, the search for
concrete versions of Ribe's theorem has fueled some of the field's
most important achievements.

The notions of type and cotype of Banach spaces are the basis of a
deep and rich theory which encompasses diverse aspects of the
local theory of Banach spaces. We refer
to~\cite{MS86}, \cite{Pisier86}, \cite{Pisier86-book}, 
\cite{T-J89}, \cite{Pisier89}, \cite{LT91}, \cite{DJT95}, \cite{Woj96}, \cite{Maurey03}
and the references therein for background on these topics. A
Banach space $X$ is said to have (Rademacher) type $p> 0$ if there
exists a constant $T<\infty$ such that for every $n$ and every
$x_1,\ldots,x_n\in X$,
\begin{eqnarray}\label{eq:def type}
\E_\e\Le\|\sum_{j=1}^n \e_j x_j\Ri\|_X^p\le T^p\sum_{j=1}^n
\|x_j\|_X^p.
\end{eqnarray}
where the expectation $\E_\e$ is with respect to a uniform choice
of signs $\e=(\e_1,\ldots,\e_n)\in \{-1,1\}^n$. $X$ is said to have
(Rademacher) cotype $q>0$ if there exists a constant $C<\infty$
such that for every $n$ and every $x_1,\ldots,x_n\in X$,
\begin{eqnarray}\label{eq:def Rademacher cotype}
\E_\e\Le\|\sum_{j=1}^n \e_j x_j\Ri\|_X^q\ge
\frac{1}{C^q}\sum_{j=1}^n \|x_j\|_X^q.
\end{eqnarray}

These notions are clearly {\em linear} notions, since their
definition involves addition and multiplication by scalars. Ribe's
theorem implies that these notions are preserved under uniform
homeomorphisms of Banach spaces, and therefore it would be
desirable to reformulate them using only distances between points
in the given Banach space. Once this is achieved, one could define
the notion of type and cotype of a metric space, and then
hopefully transfer some of the deep theory of type and cotype to
the context of arbitrary metric spaces. The need for such a theory
has recently received renewed impetus due to the discovery of
striking applications of metric geometry to theoretical computer
science (see~\cite{Mat01}, \cite{Ind01}, \cite{Linial02} and the references
therein for part of the recent developments in this direction).

Enflo's pioneering work~\cite{Enflo69}, \cite{Enflo69-groups}, \cite{Enflo70}, \cite{Enflo78} resulted in
the formulation of a nonlinear notion of type, known today as {\em Enflo type}. The basic idea is
that given a Banach space $X$ and $x_1,\ldots,x_n\in X$, one can
consider the {\em linear} function $f:\{-1,1\}^n\to X$ given by
$f(\e)=\sum_{j=1}^n\e_j x_j$. Then~\eqref{eq:def type} becomes
\begin{multline}\label{eq:enflo type norm}
\E_\e \|f(\e)-f(-\e)\|_X^p 
 \le  T^p\sum_{j=1}^n
\E_\e\Big\|f(\e_1,\ldots,\e_{j-1},\e_j,\e_{j+1},\ldots,\e_n)\\
 -f(\e_1,\ldots,\e_{j-1},-\e_j,\e_{j+1},\ldots,\e_n)\Big\|_X^p. 
\end{multline}
One can thus say that a metric space $(\M,d_\M)$ has Enflo type
$p$ if there exists a constant $T$ such that for every $n\in
\mathbb N$ and {\em every} $f:\{-1,1\}^n\to \M$,
\begin{multline}\label{eq:enflo type}
\E_\e d_\M\left(f(\e),f(-\e)\right)^p\le T^p\sum_{j=1}^n \E_\e
d_\M\Big(f(\e_1,\ldots,\e_{j-1},\e_j,\e_{j+1},\ldots,\e_n),\\
f(\e_1,\ldots,\e_{j-1},-\e_j,\e_{j+1},\ldots,\e_n)\Big)^p.
\end{multline}
There are two natural concerns about this definition. First of all,
while in the category of Banach spaces~\eqref{eq:enflo type} is
clearly a strengthening of~\eqref{eq:enflo type norm} (as we are not
restricting only to linear functions $f$), it isn't clear
whether~\eqref{eq:enflo type} follows from~\eqref{eq:enflo type
norm}. Indeed, this problem was posed  by Enflo in~\cite{Enflo78},
and in full generality it remains open. Secondly, we do not know
if~\eqref{eq:enflo type} is a useful notion, in the sense that it
yields metric variants of certain theorems from the linear theory of
type (it should be remarked here that Enflo found striking
applications of his notion of type to Hilbert's fifth problem in
infinite dimensions~\cite{Enflo69-groups}, \cite{Enflo70}, \cite{Enflo78}, and to
the uniform classification of $L_p$ spaces~\cite{Enflo69}). As we
will presently see, in a certain sense both of these issues turned
out not to be problematic. Variants of Enflo type were studied by
Gromov~\cite{Gro83} and Bourgain, Milman and Wolfson~\cite{BMW86}.
Following  \cite{BMW86} we shall say that
a metric space $(\M,d_\M)$ has BMW type $p>0$ if there exists a
constant $K<\infty$ such that for every $n\in \mathbb N$ and every
$f:\{-1,1\}^n\to \M$,
\begin{multline}\label{eq:BMW}
\E_\e d_\M(f(\e),f(-\e))^2\le K^2n^{\frac{2}{p}-1}\sum_{j=1}^n
\E_\e
d_\M\Big(f(\e_1,\ldots,\e_{j-1},\e_j,\e_{j+1},\ldots,\e_n),
\\ f(\e_1,\ldots,\e_{j-1},-\e_j,\e_{j+1},\ldots,\e_n)\Big)^2.
\end{multline}
Bourgain, Milman and Wolfson proved in~\cite{BMW86} that if a
Banach space has BMW type $p>0$ then it also has Rademacher type
$p'$ for all $0<p'<p$. They also obtained a nonlinear version of
the Maurey-Pisier theorem for type~\cite{Pisier74}, \cite{MP76}, yielding
a characterization of metric spaces which contain bi-Lipschitz
copies of the Hamming cube. In~\cite{Pisier86} Pisier proved that
for Banach spaces, Rademacher type $p$ implies Enflo type $p'$ for
every $0<p'<p$. Variants of these problems were studied by Naor
and Schechtman in~\cite{NS02}. A stronger notion of nonlinear
type, known as Markov type, was introduced by Ball~\cite{Ball92}
in his study of the {\em Lipschitz extension problem}. This
important notion has since found applications to various
fundamental problems in metric
geometry~\cite{Naor01}, \cite{LMN02}, \cite{BLMN05}, \cite{NPSS04}, \cite{MN05-proc}

Despite the vast amount of research on nonlinear type, a
nonlinear notion of cotype remained elusive. Indeed, the problem
of finding a notion of cotype which makes sense for arbitrary
metric spaces, and which coincides (or almost coincides) with the
notion of Rademacher type when restricted to Banach spaces, became
a central open problem in the field.

 There are several difficulties involved in defining nonlinear
cotype. First of all, one cannot simply reverse
inequalities~\eqref{eq:enflo type} and~\eqref{eq:BMW}, since the
resulting condition fails to hold true even for Hilbert space
(with $p=2$). Secondly, if Hilbert space satisfies an inequality
such as~\eqref{eq:enflo type}, then it must satisfy the same
inequality where the distances are raised to any power $0<r<p$.
This is because Hilbert space, equipped with the metric
$\|x-y\|^{r/p}$, is isometric to a subset of Hilbert space
(see~\cite{Schoenberg38}, \cite{WW75}). In the context of nonlinear type,
this observation makes perfect sense, since if a Banach space has
type $p$ then it also has type $r$ for every $0<r<p$. But, this is
no longer true for cotype (in particular, no Banach space has
cotype less than $2$). One viable definition of cotype of a metric
space $X$ that was suggested in the early 1980s is the following:
Let $\M$ be a metric space, and denote by $\Lip(\M)$ the Banach
space of all real-valued Lipschitz functions on $\M$, equipped
with the Lipschitz norm. One can then define the nonlinear cotype
of $\M$ as the (Rademacher) cotype of the (linear) dual
$\Lip(\M)^*$. This is a natural definition when $\M$ is a Banach
space, since we can view $\Lip(\M)$ as a nonlinear substitute for
the dual space $\M^*$ (note that in~\cite{Lin64} it is shown that
there is a norm $1$ projection from $\Lip(\M)$ onto $\M^*$). With
this point of view, the above definition of cotype is natural due
to the principle of local reflexivity~\cite{LR69}, \cite{JRZ71}.
Unfortunately, Bourgain~\cite{Bourgain86-trees} has shown that
under this definition subsets of $L_1$ need not have finite
nonlinear cotype (while $L_1$ has cotype $2$). Additionally, the
space $\Lip(M)^*$ is very hard to compute, for example it is an
intriguing open problem whether even the unit square $[0,1]^2$ has
nonlinear cotype $2$ under the above definition.

In this paper we introduce a notion of cotype of metric spaces,
and show that it coincides with Rademacher cotype when restricted
to the category of Banach spaces. Namely, we introduce the
following concept:

\begin{definition}[Metric cotype]\label{def:cotype} Let $(\M,d_\M)$ be a
 metric space and\break $q>0$. The space
$(\M,d_\M)$ is said to have {\it metric cotype $q$ with constant} $\Gamma$ if for
every integer $n\in \mathbb N$, there exists an even integer $m$,
such that for every $f:\Z_m^n\to \M$,
\begin{eqnarray}\label{eq:def cotype}
\sum_{j=1}^n\E_x\left[d_\M\left(f\left(x+\frac{m}{2}e_j\right),f(x)\right)^q\right]\le
\Gamma^q m^q\E_{\e,x}\left[d_\M(f(x+\e),f(x))^q\right],
\end{eqnarray}
where the expectations above are taken with respect to uniformly
chosen $x\in \Z_m^n$ and $\e\in\{-1,0,1\}^n$ (here, and in what
follows we denote by $\{e_j\}_{j=1}^n$ the standard basis of
$\R^n$). The smallest constant $\Gamma$ with which
inequality~\eqref{eq:def cotype} holds true is denoted
$\Gamma_q(\M)$.
\end{definition}

Several remarks on Definition~\ref{def:cotype} are in order. First
of all, in the case of Banach spaces, if we apply
inequality~\eqref{eq:def cotype} to linear functions
$f(x)=\sum_{j=1}^n x_j v_j$, then by homogeneity $m$ would cancel,
and the resulting inequality will simply become the Rademacher
cotype $q$ condition (this statement is not precise due to the
fact that addition on $\Z_m^n$ is performed modulo $m$ --- see
Section~\ref{section:easy direction} for the full argument). Secondly,
it is easy to see that in any metric space which contains at least
two points, inequality~(\ref{eq:def cotype}) forces the scaling
factor $m$ to be large (see Lemma~\ref{lem:lower m}) --- this is an
essential difference between Enflo type and metric cotype.
Finally, the averaging over $\e\in \{-1,0,1\}^n$ is natural here,
since this forces the right-hand side of~\eqref{eq:def cotype}  to be a
uniform average over all pairs in $\Z_m^n$ whose distance is at
most $1$ in the $\ell_\infty$ metric.
\medbreak
 
The following theorem is the main result of this paper:

\begin{theorem}\label{thm:cotype} Let $X$ be a Banach space{\rm ,} and $q\in [2,\infty)$. Then
$X$ has metric cotype $q$ if and only if $X$ has Rademacher cotype
$q$. Moreover{\rm ,}
$$
\frac{1}{2\pi}C_q(X)\le \Gamma_q(X)\le 90C_q(X).
$$
\end{theorem}

Apart from settling the nonlinear cotype problem described above,
this notion has various applications. Thus, in the remainder of
this paper we proceed to study metric cotype and some of its
applications, which we describe below. We believe that additional
applications of this notion and its variants will be discovered in
the future. In particular, it seems worthwhile to study the
interaction between metric type and metric cotype (such as in
Kwapien's theorem~\cite{Kwapien72}), the possible ``Markov"
variants of metric cotype (\`a la Ball~\cite{Ball92}) and their
relation to the Lipschitz extension problem, and the relation
between metric cotype and the nonlinear Dvoretzky theorem
(see~\cite{BFM86}, \cite{BLMN05} for information about the nonlinear
Dvoretzky theorem, and~\cite{FLM77} for the connection between
cotype and Dvoretzky's theorem).

\subsection{Some applications of metric cotype}
\smallbreak

  \subsubsection*{1) \textsl{A nonlinear version of the Maurey-Pisier
theorem}.} Given two metric spaces $(\M,d_\M)$ and $(\N,d_\N)$,
and an injective mapping $f:\M\hookrightarrow \N$, we denote the
{\em distortion} of $f$ by
$$
\dist(f):= \|f\|_{\Lip}\cdot\|f^{-1}\|_{\Lip}=\sup_{\substack{x,y\in
\M\\ x\neq y}} \frac{d_\N(f(x),f(y))}{d_\M(x,y)}\cdot
\sup_{\substack{x,y\in \M\\ x\neq y}}
\frac{d_\M(x,y)}{d_\N(f(x),f(y))}.
$$
The smallest distortion with which $\M$ can be embedded into $\N$
is denoted $c_\N(\M)$; i.e.,
$$
c_\N(\M):= \inf\{\dist(f):\ f:\M\hookrightarrow \N\}.
$$
If $c_\N(\M)\le \alpha$ then we sometimes use the notation $ \M
\overset{\alpha}{\hookrightarrow} \N$. When $\N=L_p$ for some
$p\ge 1$, we write $c_\N(\cdot)=c_p(\cdot)$.

For a Banach space $X$ write
$$
p_X=\sup\{p\ge 1:\ T_p(X)<\infty\}\quad \mathrm{and}\quad
q_X=\inf\{q\ge 2:\ C_q(X)<\infty\}.
$$
$X$ is said to have nontrivial type if $p_X>1$, and $X$ is said
to have nontrivial cotype if $q_X<\infty$.

In~\cite{Pisier74} Pisier proved that $X$ has no nontrivial type
if and only if for every $n\in \mathbb N$ and every $\e>0$,
$\ell_1^n\overset{1+\e}{\hookrightarrow} X$. A nonlinear analog
of this result was proved by Bourgain, Milman and
Wolfson~\cite{BMW86} (see also Pisier's proof in~\cite{Pisier86}).
They showed that a metric space $\M$ does not have BMW type larger
than $1$ if and only if for every $n\in \mathbb N$ and every
$\e>0$,
$(\{0,1\}^n,\|\cdot\|_1)\overset{1+\e}{\hookrightarrow}\M$.
In~\cite{MP76} Maurey and Pisier proved that a Banach space $X$
has no nontrivial cotype if and only for every $n\in \mathbb N$
and every $\e>0$, $\ell_\infty^n \overset{1+\e}{\hookrightarrow}
X$. To obtain a nonlinear analog of this theorem we need to
introduce a variant of metric cotype (which is analogous to the
variant of Enflo type that was used in [11].

\begin{definition}[Variants of metric cotype \`a la Bourgain, Milman and\break Wolfson]
   Let $(\M,d_\M)$ be a metric space and
$1\le p\le q$. We denote by $\Gamma^{(p)}_q(\M)$ the least constant
$\Gamma$ such that for every integer $n\in \mathbb N$ there exists
an even integer $m$, such that for every $f:\Z_m^n\to \M$,
\begin{multline}\label{eq:def weak cotype}
\sum_{j=1}^n\E_x\left[d_\M\left(f\left(x+\frac{m}{2}e_j\right),f(x)\right)^p\right]\\
\le
\Gamma^p m^p
n^{1-\frac{p}{q}}\E_{\e,x}\left[d_\M(f(x+\e),f(x))^p\right],
\end{multline}
where the expectations above are taken with respect to uniformly
chosen $x\in \Z_m^n$ and $\e\in\{-1,0,1\}^n$. 
Note that
$\Gamma^{(q)}_q(\M)=\Gamma_q(\M)$. When $1\le p<q$ we shall refer
to~\eqref{eq:def weak cotype} as a weak metric cotype $q$
inequality with exponent $p$ and constant $\Gamma$.
\end{definition}
\setcounter{theorem}{3}

The following theorem is analogous to Theorem~\ref{thm:cotype}.

\begin{theorem}\label{thm:weak cotype}
Let $X$ be a Banach space{\rm ,} and assume that for some $1\le p< q${\rm ,}
$\Gamma_q^{(p)}(X)<\infty$. Then $X$ has cotype $q'$ for every
$q'>q$. If $q=2$ then $X$ has cotype $2$. On the other hand{\rm ,}
$$
\Gamma_q^{(p)}(X)\le c_{pq}C_q(X),
$$
where $c_{pq}$ is a universal constant depending only on $p$ and
$q$.
\end{theorem}

In what follows, for $m,n\in \mathbb N$ and $p\in [1,\infty]$ we
let $[m]_p^n$ denote the set $\{0,1,\ldots,m\}^n$, equipped with
the metric induced by $\ell_p^n$. The following theorem is a
metric version of the Maurey-Pisier theorem (for cotype):

\begin{theorem}\label{thm:MPcotype} Let $\M$ be a metric space
such that $\Gamma_q^{(2)}(\M)=\infty$ for all $q<\infty$. Then for
every $m,n\in \mathbb N$ and every $\e>0${\rm ,}
$$
[m]_\infty^n \overset{1+\e}{\hookrightarrow}\M.
$$
\end{theorem}

We remark that in~\cite{MP76} Maurey and Pisier prove a stronger
result, namely that for a Banach space $X$, for every $n\in
\mathbb N$ and every $\e>0$, $\ell_{p_X}^n
\overset{1+\e}{\hookrightarrow} X$ and $\ell_{q_X}^n
\overset{1+\e}{\hookrightarrow} X$. Even in the case of nonlinear
type, the results of Bourgain, Milman and Wolfson yield an
incomplete analog of this result in the case of BMW type greater
than $1$. The same phenomenon seems to occur when one tries to
obtain a nonlinear analog of the full Maurey-Pisier theorem for
cotype. We believe that this issue deserves more attention in
future research.

\subsubsection*{{\rm 2)} \textsl{Solution of a problem posed by  Arora{\rm ,} Lov\'{a}sz{\rm ,} Newman{\rm ,}
Rabani{\rm ,}\break Rabinovich and Vempala.}}
The following question appears in \cite[Conj.~5.1]{ALNRRV05}:

\begin{quote}
Let $\mathcal F$ be a \emph{baseline} metric class which does not contain
all finite metrics with distortion arbitrarily close to $1$. Does
this imply that there exists $\alpha>0$ and arbitrarily large
$n$-point metric spaces $\M_n$ such that for every $\N\in \mathcal
F$, $c_\N(\M_n)\ge (\log n)^\alpha$?
\end{quote}

We refer to~\cite[\S 2]{ALNRRV05} for the definition of baseline
metrics, since we will not use this notion in what follows. We also
refer to~\cite{ALNRRV05} for background and motivation from
combinatorial optimization for this problem, where several partial
results in this direction are obtained. An extended abstract of the
current paper~\cite{MN06} also contains more information on the
connection to Computer Science. Here we apply metric cotype to
settle this conjecture positively, without any restriction on the
class $\mathcal F$.

To state our result we first introduce some notation. If $\F$ is a
family of metric spaces we write
$$
c_\F(\N)=\inf \left\{c_\M(\N): \M\in \F\right\}.
$$
For an integer $n\ge 1$ we define
$$
\mathcal D_n(\F)= \sup\{c_\F(\N):\ \N\ \text{is a metric space},\
|\N|\le n\}.
$$
Observe that if, for example, $\F$ consists of all the subsets of
Hilbert space (or $L_1$), then Bourgain's embedding
theorem~\cite{Bou85} implies that $\mathcal D_n(\F)=O(\log n)$.

For $K>0$ we define the $K$-cotype (with exponent $2$) of a family
of metric spaces $\F$ as
$$
q_\F^{(2)}(K)=\sup_{\M\in \F} \inf \left\{q\in (0,\infty]:\
\Gamma_q^{(2)}(\M)\le K\right\}.
$$
Finally we let
$$
q^{(2)}_\F=\inf_{\infty>K>0} q_\F^{(2)}(K).
$$

The following theorem settles positively the problem stated above:
\begin{theorem}\label{thm:dicho}
Let $\F$ be a family of metric spaces. Then the following
conditions are equivalent\/{\rm :}\/
\begin{enumerate} \itemsep 0mm
\item 
There exists a finite metric space $\M$ for which $c_{\F}(\M)>1$.
\item  $q_\F^{(2)}<\infty$. \item There exists $0<\alpha<\infty$
such that $\mathcal D_n(\F)=\Omega\left((\log n)^\alpha\right)$.
\end{enumerate}
\end{theorem}

\subsubsection*{ 3) \textsl{A quantitative version of Matou{\hskip.75pt\rm \v{{\hskip-5.75pt\it s}}}ek\/{\rm '}\/s {\rm BD}
Ramsey theorem.}} In~\cite{Mat92} Matou\v{s}ek proved the
following result, which he calls the Bounded Distortion (BD)
Ramsey theorem. We refer to~\cite{Mat92} for motivation and
background on these types  of results.

\vskip8pt{{\sc Theorem 1.7} {\rm (Matou\v{s}ek's BD Ramsey theorem).}}
{\it Let $ X$ be a finite metric space and $\e>0${\rm ,} $\gamma>1$. Then
there exists a metric space $ Y=Y(X,\e,\gamma)${\rm ,} such that for
every metric space} $Z${\rm ,}
$$
 c_Z(Y)<\gamma\implies c_Z(X)<1+\e.
$$
\vskip6pt

We obtain a new proof of Theorem~1.7, which is
quantitative and concrete:

\vskip6pt {{\sc Theorem 1.8} {\rm (Quantitative version of Matou\v{s}ek's BD Ramsey
theorem).}} {\it There exists a universal constant $C$
with the following properties. Let $X$ be an $n$-point metric
space and $\e\in (0,1)${\rm ,} $\gamma>1$. Then for every integer $N\ge
(C\gamma)^{2^{5A}}${\rm ,} where
$$
A=\max\left\{\frac{4\diam(X)}{\e\cdot \min_{x\neq y}
d_X(x,y)},n\right\},
$$
if a metric space $Z$ satisfies $c_Z(X)>1+\e$ then{\rm ,}
$c_Z\left(\left[N^5\right]_\infty^N\right)>\gamma$.}
\vskip6pt
\setcounter{theorem}{8}

We note that Matou\v{s}ek's argument in~\cite{Mat92} uses Ramsey
theory, and is nonconstructive (at best it can yield tower-type
bounds on the size of $Z$, which are much worse than what the
cotype-based approach gives).

\subsubsection*{{\rm 4)} \textsl{Uniform embeddings and Smirnov\/{\rm '}\/s problem}.} Let
$(\M,d_\M)$ and $(\N,d_\N)$ be metric spaces. A mapping $f:\M\to
\N$ is called a {\em uniform embedding} if $f$ is injective, and
both $f$ and $f^{-1}$ are uniformly continuous. There is a large
body of work on the uniform classification of metric spaces --- we
refer to the survey article~\cite{Lin98}, the book~\cite{BL00},
and the references therein for background on this topic. In spite
of this, several fundamental questions remain open. For example,
it was not known for which values of $0<p, q<\infty$,  $L_p$ embeds
uniformly into $L_q$. As we will presently see, our results yield
a complete characterization of these values of $p,q$.

In the late 1950's Smirnov asked whether every separable metric
space embeds uniformly into $L_2$ (see~\cite{Gorin59}). Smirnov's
problem was settled negatively by Enflo in~\cite{Enflo69-smirnov}.
Following Enflo, we shall say that a metric space $\M$ is a {\em
universal uniform embedding space} if every separable metric space
embeds uniformly into $\M$. Since every separable metric space is
isometric to a subset of $C[0,1]$, this is equivalent to asking
whether $C[0,1]$ is uniformly homeomorphic to a subset of $\M$
(the space $C[0,1]$ can be replaced here by $c_0$ due to Aharoni's
theorem~\cite{Aha74}). Enflo proved that $c_0$ does not uniformly
embed into Hilbert space. In~\cite{AMM85}, Aharoni, Maurey and
Mityagin systematically studied metric spaces which are uniformly
homeomorphic to a subset of Hilbert space, and obtained an elegant
characterization of Banach spaces which are uniformly homeomorphic
to a subset of $L_2$. In particular, the results of~\cite{AMM85}
imply that for $p>2$, $L_p$ is not uniformly homeomorphic to a
subset of $L_2$.

Here we prove that in the class of Banach spaces with nontrivial
type, if $Y$ embeds uniformly into $X$, then $Y$ inherits the
cotype of $X$. More precisely:

\begin{theorem}\label{thm:uniform} Let $X$ be a
Banach space with nontrivial type. Assume that $Y$ is a Banach
space which uniformly embeds into $X$. Then $q_Y\le q_X$.
\end{theorem}

As a corollary, we complete the characterization of the values of
$0<p$,\break $q<\infty$ for which $L_p$ embeds uniformly into $L_q$:

\begin{theorem}\label{thm:uniformL_p}
For $p,q> 0${\rm ,} $L_p$ embeds uniformly into $L_q$ if and only if
$p\le q$ or $q\le p\le 2$.
\end{theorem}

We believe that the assumption that $X$ has nontrivial type in
Theorem~\ref{thm:uniform} can be removed --- in
Section~\ref{section:problems} we present a concrete problem which
would imply this fact. If true, this would imply that cotype is
preserved under uniform embeddings of Banach spaces. In
particular, it would follow that a universal uniform embedding
space cannot have nontrivial cotype, and thus by the
Maurey-Pisier theorem~\cite{MP76} it must contain
$\ell_\infty^n$'s with distortion uniformly bounded in $n$.

\subsubsection*{{\rm 5)} \textsl{Coarse embeddings}.} Let $(\M,d_\M)$ and
$(\N,d_\N)$ be metric spaces. A mapping $f:\M\to \N$ is called a
{\em coarse embedding} if there exists two nondecreasing
functions $\alpha,\beta:[0,\infty)\to[0,\infty)$ such that
$\lim_{t\to\infty} \alpha(t)=\infty$, and for every $x,y\in \M$,
$$
\alpha(d_\M(x,y))\le d_\N(f(x),f(y))\le \beta(d_\M(x,y)).
$$
This (seemingly weak) notion of embedding was introduced by Gromov
(see \cite{Gro99}), and has  several important geometric
applications. In particular, Yu~\cite{Yu00} obtained a striking
connection between the Novikov and Baum-Connes conjectures and
coarse embeddings into Hilbert spaces. In~\cite{KY04} Kasparov and
Yu generalized this to coarse embeddings into arbitrary uniformly
convex Banach spaces. It was unclear, however, whether this is
indeed a strict generalization, i.e. whether or not the existence
of a coarse embedding into a uniformly convex Banach space implies
the existence of a coarse embedding into a Hilbert space. This was
resolved by Johnson and Randrianarivony in~\cite{JR04}, who proved
that for $p>2$, $L_p$ does not coarsely embed into $L_2$.
In~\cite{Ran04}, Randrianarivony proceeded to obtain a
characterization of Banach spaces which embed coarsely into $L_2$,
in the spirit of the result of Aharoni, Maurey and
Mityagin~\cite{AMM85}. There are very few known methods of proving
coarse nonembeddability results. Apart from the
papers~\cite{JR04}, \cite{Ran04} quoted above, we refer
to~\cite{Gro03}, \cite{DGLY02}, \cite{Oza04} for results of this type. Here we use
metric cotype to prove the following coarse variants of
Theorem~\ref{thm:uniform} and Theorem~\ref{thm:uniformL_p}, which
generalize, in particular, the theorem of Johnson and
Randrianarivony.

\begin{theorem}\label{thm:coarse}
Let $X$ be a Banach space with nontrivial type. Assume that $Y$
is a Banach space which coarsely embeds into $X$. Then $q_Y\le
q_X$. In particular{\rm ,} for $p,q> 0${\rm ,} $L_p$ embeds coarsely into
$L_q$ if and only if $p\le q$ or $q\le p\le 2$.
\end{theorem}

 \subsubsection*{\textrm{6)}  \textsl{Bi-Lipschitz embeddings of the integer lattice.}}
Bi-Lipschitz embeddings of the integer lattice $[m]_p^n$ were
investigated by Bourgain in~\cite{Bourgain87} and by the present
authors in~\cite{MN05-proc} where  it was shown that
if $2\le p<\infty$ and $Y$ is a Banach space which admits an
equivalent norm whose modulus of uniform convexity has power type
$2$, then
\begin{equation}\label{eq:phase}
c_Y\left([m]_p^n\right)=\Theta\left(\min\left\{n^{\frac12-\frac{1}{p}},m^{1-\frac{2}{p}}\right\}\right).
\end{equation}
The implied constants in the above asymptotic equivalence
depend on $p$ and on the $2$-convexity constant of $Y$. Moreover,
it was shown in~\cite{MN05-proc} that
$$
c_Y([m]_\infty^n)=\Omega\left(\min\left\{\sqrt{\frac{n}{\log
n}},\frac{m}{\sqrt{\log m}}\right\}\right).
$$
It was conjectured in~\cite{MN05-proc} that the logarithmic terms
above are unnecessary. Using our results on metric cotype we
settle this conjecture positively, by proving the following
general theorem:

\begin{theorem}\label{thm:infty grid} Let $Y$ be a Banach space with nontrivial type which has
cotype $q$. Then
$$
c_Y([m]_\infty^n)=\Omega\left(\min\left\{n^{1/q},m\right\}\right).
$$
\end{theorem}

Similarly, our methods imply that~\eqref{eq:phase} holds true for
any Banach space $Y$ with nontrivial type and cotype $2$ (note
that these conditions are strictly weaker than being $2$-convex,
as shown e.g. in~\cite{LTII77}). Moreover, it is possible to
generalize the lower bound in~\eqref{eq:phase} to Banach spaces
with nontrivial type, and cotype $2\le q\le p$, in which case the
lower bound becomes
$\min\left\{n^{\frac{1}{q}-\frac{1}{p}},m^{1-\frac{q}{p}}\right\}$.

\subsubsection*{{\rm 7)} \textsl{Quadratic inequalities on the cut-cone.}} An
intriguing aspect of Theorem~\ref{thm:cotype} is that $L_1$ has
metric cotype $2$. Thus, we obtain a nontrivial inequality on
$L_1$ which involves distances {\em squared}. To the best of our
knowledge, all the known nonembeddability results for $L_1$ are
based on Poincar\'e type inequalities in which distances are
raised to the power $1$. Clearly, any such inequality reduces to
an inequality on the real line. Equivalently, by the cut-cone
representation of $L_1$ metrics (see~\cite{DL97}) it is enough to
prove any such inequality for {\em cut metrics}, which are
particularly simple. Theorem~\ref{thm:cotype} seems to be the
first truly ``infinite dimensional" metric inequality in $L_1$, in
the sense that its nonlinearity does not allow a straightforward
reduction to the one-dimensional case. We believe that
understanding such inequalities on $L_1$ deserves further
scrutiny, especially as they hint at certain nontrivial (and
nonlinear) interactions between cuts.

\section{Preliminaries and notation}

We start by setting   notation and conventions. Consider the
standard $\ell_\infty$ Cayley graph on $\Z_m^n$, namely $x,y\in
\Z_m^n$ are joined by an edge if and only if they are distinct and
$x-y\in \{-1,0,1\}^n$. This induces a shortest-path metric on
$\Z_m^n$ which we denote by $d_{\Z_m^n}(\cdot,\cdot)$.
Equivalently, the metric space $(\Z_m^n,d_{\Z_m^n})$ is precisely
the quotient $(\Z^n,\|\cdot\|_\infty)/(m\Z)^n$ (for background on
quotient metrics see~\cite{BH99}, \cite{Gro99}). The ball of radius $r$
around $x\in \Z_m^n$ will be denoted $B_{\Z_m^n}(x,r)$. We denote
by $\mu$ the normalized counting measure on $\Z_m^n$ (which is
clearly the Haar measure on this group). We also denote by
$\sigma$ the normalized counting measure on $\{-1,0,1\}^n$. In
what follows, whenever we average over uniformly chosen signs
$\e\in\{-1,1\}^n$ we use the probabilistic notation $\E_\e$ (in
this sense we break from the notation used in the introduction,
for the sake of clarity of the ensuing arguments).

In what follows all Banach spaces are assumed to be over the
complex numbers $\mathbb C$. All of our results hold for real
Banach spaces as well, by a straightforward complexification
argument.

Given a Banach space $X$ and $p,q\in [1,\infty)$ we denote by
$C_q^{(p)}(X)$ the infimum over all constants $C>0$ such that for
every integer $n\in \mathbb N$ and every $x_1,\ldots,x_n\in X$,
\begin{eqnarray}\label{eq:pass to p}
\Biggl(\E_\e\Biggl\|\sum_{j=1}^n \e_j
x_j\Biggr\|_X^p\Biggr)^{1/p}\ge \frac{1}{C}\Biggl(\sum_{j=1}^n
\|x_j\|_X^q\Biggr)^{1/q}.
\end{eqnarray}
Thus, by our previous notation, $C_q^{(q)}(X)=C_q(X)$. Kahane's
inequality~\cite{kahane64} says that for $1\le p,q<\infty$ there
exists a constant $1\le A_{pq}<\infty$ such that for every Banach
space $X$, every integer $n\in \mathbb N$, and every
$x_1,\ldots,x_n\in X$,
\begin{eqnarray}\label{eq:kahane}
 \Biggl(\E_\e\Biggl\|\sum_{j=1}^n
\e_j x_j\Biggr\|_X^p\Biggr)^{1/p}\le
A_{pq}\Biggl(\E_\e\Biggl\|\sum_{j=1}^n \e_j
x_j\Biggr\|_X^q\Biggr)^{1/q}.
\end{eqnarray}
Where clearly $A_{pq}=1$ if $p\le q$, and for every $1\le
q<p<\infty$, $A_{pq}=O\left(\sqrt{p}\right)$ (see~\cite{Tal88}).
It follows in particular from~\eqref{eq:kahane} that if $X$ has
cotype $q$ then for every $p\in [1,\infty)$,
$C_q^{(p)}(X)=O_{p,q}(C_q(X))$, where the implied constant may
depend on $p$ and $q$.

Given $A\!\subseteq\!\{1,\ldots,n\}$,
we consider the Walsh functions \hbox{$W_A:\{-1,1\}^n \to \mathbb{C}$,} defined as
\[ W_A(\e_1,\ldots,\e_m)=\prod_{j\in A} \e_j .\]
Every $f:\{-1,1\}^n\to X$ can be written as
$$
f(\e_1,\ldots,\e_n)=\sum_{A\subseteq \{1,\ldots,n\}}\widehat
f(A)W_A(\e_1,\ldots,\e_n),
$$
where $\widehat f(A)\in X$ are given by
$$
\widehat f(A)=\E_\e \Bigl(f(\e)W_A(\e) \Bigr).
$$
The {\em Rademacher projection} of $f$ is defined by
$$
\Rad(f)=\sum_{j=1}^n \widehat f(A)W_{\{j\}}.
$$
The $K$-convexity constant of $X$, denoted $K(X)$, is the smallest
constant $K$ such that for every $n$ and every $f:\{-1,1\}^n\to
X$,
$$
\E_\e\|\Rad(f)(\e)\|_X^2\le K^2 \E_\e \|f(\e)\|_X^2.
$$
In other words,
$$
K(X)=\sup_{n\in \mathbb N} \|\Rad\|_{L_2(\{-1,1\}^n,X)\to
L_2(\{-1,1\}^n,X)}.
$$
$X$ is said to be $K$-convex if $K(X)<\infty$. More generally, for
$p\ge 1$ we define
$$
K_p(X)=\sup_{n\in \mathbb N} \|\Rad\|_{L_p(\{-1,1\}^n,X)\to
L_p(\{-1,1\}^n,X)}.
$$
It is a well known consequence of Kahane's inequality and duality
that for every $p>1$,
$$
K_p(X)\le O\Biggl(\frac{p}{\sqrt{p-1}}\Biggr)\cdot K(X).
$$

The following deep theorem was proved by Pisier in~\cite{Pis82}:

\begin{theorem}[Pisier's $K$-convexity theorem~\cite{Pis82}]
Let $X$ be a Banach space. Then
$$
q_X>1\iff K(X)<\infty.
$$
\end{theorem}

Next, we recall some facts concerning Fourier analysis on the
group $\Z_m^n$. Given $k=(k_1,\ldots,k_n)\in \Z_m^n$ we consider
the Walsh function $W_k:\Z_m^n\to \mathbb C$:
$$
W_{k}(x)=\exp\Biggl(\frac{2\pi i}{m}\sum_{j=1}^m k_jx_j\Biggr).
$$
Then, for any Banach space $X$, any $f:Z_m^n\to X$ can be
decomposed as follows:
$$
f(x)=\sum_{k\in \mathbb Z_m^n} W_k(x)\widehat f(k),
$$
where
$$
\widehat f(k)=\int_{\Z_m^n} f(y)\overline{W_k(y)}d\mu(y)\in X.
$$
If $X$ is a Hilbert space then Parseval's identity becomes:
$$
\int_{Z_m^n} \|f(x)\|_X^2d\mu(x)=\sum_{k\in \Z_m^n}
\left\|\widehat f(k)\right\|_X^2.
$$

\subsection{Definitions and basic facts related to metric cotype}

\begin{definition}
Given $1\le p\le q$, an integer $n$ and an even integer $m$, let
$\Gamma_q^{(p)}(\M;n,m)$ be the infimum over all $\Gamma>0$ such
that for every $f:\Z_m^n\to \M$, \begin{multline}\label{eq:two
parameter} \sum_{j=1}^n \int_{\Z_m^n}
d_\M\left(f\left(x+\frac{m}{2}e_j\right),f(x)\right)^pd\mu(x)\\
\le
\Gamma^pm^pn^{1-\frac{p}{q}}\int_{\{-1,0,1\}^n}\int_{\Z_m^n}d_\M\left(f\left(x+\e\right),f(x)\right)^pd\mu(x)d\sigma(\e).
\end{multline}
When $p=q$ we write $\Gamma_q(\M; n,m):= \Gamma_q^{(q)} (\M; n,m)$ .
With this notation,
$$
\Gamma_q^{(p)}(\M)=\sup_{n\in \mathbb N}\inf_{m\in 2\mathbb
N}\Gamma_q^{(p)}(\M;n,m).
$$

We also denote by $m_q^{(p)}(\M;n,\Gamma)$ the smallest even
integer $m$ for which~\eqref{eq:two parameter} holds. As usual,
when $p=q$ we write $m_q(\M;n,\Gamma):= m_q^{(q)}(\M;n,\Gamma)$.
\end{definition}

The following lemma shows that for nontrivial metric spaces $\M$,\break
$m_q(\M;n,\Gamma)$ must be large.

\begin{lemma}\label{lem:lower m}
Let $(\M,d_\M)$ be a metric space which contains at least two
points. Then for every integer $n${\rm ,} every $\Gamma>0${\rm ,} and every
$p,q>0${\rm ,}
$$
m_q^{(p)}(\M;n,\Gamma)\ge \frac{n^{1/q}}{\Gamma}.
$$
\end{lemma}

\begin{proof}  Fix $u,v\in \M$, $u \ne v$, and without loss of generality
normalize the metric so that $d_\M(u,v)=1$. Denote
$m=m_q^{(p)}(\M;n,\Gamma)$. Let $f:\Z_m^n\to \M$ be the random
mapping such that for every $x\in \Z_m^n$,
$\Pr[f(x)=u]=\Pr[f(x)=v]=\frac12$, and $\{f(x)\}_{x\in \Z_m^n}$
are independent random variables. Then for every distinct $x,y\in
\Z_m^n$, $\E \left[d_\M(f(x),f(y))^p\right]=\frac12$. Thus, the
required result follows by applying~\eqref{eq:two parameter} to $f$  and taking expectation.
\end{proof}

\begin{lemma}\label{lem:multip} For every two integers $n,k${\rm ,} and every even integer
$m${\rm ,}
$$
\Gamma_q^{(p)}(\M;n,km)\le \Gamma_q^{(p)}(\M;n,m).
$$
\end{lemma}

\begin{proof}  Fix $f:\Z_{km}^n\to \M$. For every $y\in \Z_k^n$
define $f_y: \Z_m^n\to \M$ by
$$
f_y(x)=f(kx+y).
$$
Fix $\Gamma> \Gamma_q^{(p)}(\M;n,m)$. Applying the definition of
$\Gamma_q^{(p)}(\M;n,m)$ to $f_y$, we get that
\begin{multline*}
\sum_{j=1}^n \int_{\Z_m^n}
d_\M\left(f\left(kx+\frac{km}{2}e_j+y\right),f(kx+y)\right)^pd\mu_{\Z_m^n}(x)\\\le
\Gamma^pm^pn^{1-\frac{p}{q}}\int_{\{-1,0,1\}^n}\int_{\Z_m^n}d_\M\left(f\left(kx+k\e+y\right),f(kx+y)\right)^pd\mu_{\Z_m^n}(x)d\sigma(\e).
\end{multline*}
Integrating this inequality with respect to $y\in \Z_k^n$ we see
that

\begin{small}
\begin{eqnarray*}
&& \sum_{j=1}^n \int_{\Z_{km}^n}
d_\M\left(f\left(z+\frac{km}{2}e_j\right),f(z)\right)^pd\mu_{\Z_{km}^n}(z)\\&=&\sum_{j=1}^n
\int_{\Z_k^n}\int_{\Z_m^n}
d_\M\left(f\left(kx+\frac{km}{2}e_j+y\right),f(kx+y)\right)^pd\mu_{\Z_m^n}(x)d\mu_{\Z_k^n}(y)\\
&\le&
\Gamma^pm^pn^{1-\frac{p}{q}}\hskip-4pt\int_{\{-1,0,1\}^n}\int_{\Z_k^n}\int_{\Z_m^n}d_\M\left(f\left(kx+k\e+y\right),f(kx+y)\right)^pd\mu_{\Z_m^n}(x)d\mu_{\Z_k^n}(y)d\sigma(\e)\\
&=&
\Gamma^pm^pn^{1-\frac{p}{q}}\hskip-4pt\int_{\{-1,0,1\}^n}\int_{\Z_{km}^n}d_\M\left(f\left(z+k\e\right),f(z)\right)^pd\mu_{\Z_{km}^n}(z)d\sigma(\e)\\
&\le&
\Gamma^pm^pn^{1-\frac{p}{q}}\hskip-4pt\int_{\{-1,0,1\}^n}\int_{\Z_{km}^n}k^{p-1}\sum_{s=1}^k
d_\M\left(f\left(z+s\e\right),f(z+(s-1)\e)\right)^pd\mu_{\Z_{km}^n}(z)d\sigma(\e)\\
&=&\Gamma^p(km)^pn^{1-\frac{p}{q}}\hskip-4pt\int_{\{-1,0,1\}^n}\int_{\Z_{km}^n}
d_\M\left(f\left(z+\e\right),f(z)\right)^pd\mu_{\Z_{km}^n}(z)d\sigma(\e). \qquad \qquad \qedhere
\end{eqnarray*}\end{small}
\end{proof}

\begin{lemma}\label{lem:monotone} Let $k,n$ be integers such that $k\le n${\rm ,} and let
$m$ be an even integer. Then
\[
\Gamma_q^{(p)}(\M;k,m)\le
\left(\frac{n}{k}\right)^{1-\frac{p}{q}}\cdot\Gamma_q^{(p)}(\M;n,m).
\]
\end{lemma}
\begin{proof}
  Given an 
$f:\Z_m^k\to \M$, we define an $\M$-valued function on
$\Z_m^n\cong \Z_m^k\times \Z_m^{n-k}$ by $g(x,y)=f(x)$. Applying
the definition $\Gamma_q^{(p)}(\M;n,m)$ to $g$ yields the required
inequality.
\end{proof}

We end this section by recording some general inequalities which
will be used in the ensuing arguments. In what follows $(\M,d_\M)$
is an arbitrary metric space.

\begin{lemma}\label{lem:pass to diagonals}
For every $f:\Z_m^n\to \M${\rm ,}
\begin{multline*}
\sum_{j=1}^n\int_{\Z_m^n}
d_\M\left(f(x+e_j),f(x)\right)^pd\mu(x)\\
\le 3\cdot2^{p-1}n\cdot
\int_{\{-1,0,1\}^n}\int_{\Z_m^n}
d_\M\left(f(x+\e),f(x)\right)^pd\mu(x)d\sigma(\e).
\end{multline*}
\end{lemma}

\begin{proof}  For every $x\in \Z_m^n$ and $\e\in \{-1,0,1\}^n$,
$$
d_\M(f(x+e_j),f(x))^p\le 2^{p-1}
d_\M(f(x+e_j),f(x+\e))^p+2^{p-1}d_\M(f(x+\e),f(x))^p.
$$
Thus
\begin{align*}
 \frac23\int_{\Z_m^n} &
d_\M\left(f(x+e_j),f(x)\right)^pd\mu(x)\\ 
 &= \sigma(\{\e\in
\{-1,0,1\}^n:\ \e_j\neq -1\})
\cdot \int_{\Z_m^n}
d_\M\left(f(x+e_j),f(x)\right)^pd\mu(x)\\
&\le2^{p-1}\int_{\{\e\in
\{- 1,0,1\}^n:\ \e_j\neq
-1\}}\int_{\Z_m^n}\Big(d_\M\left(f(x+e_j),f(x+\e)\right)^p\\
&\qquad +d_\M(f(x+\e),f(x))^p\Big)d\mu(x)d\sigma(\e)\\
& =2^{p-1}\int_{\{\e\in
\{- 1,0,1\}^n:\ \e_j\neq
1\}}\int_{\Z_m^n}d_\M(f(y+\e),f(y))^pd\mu(y)d\sigma(\e)\\
&\qquad +2^{p-1}\int_{\{\e\in
\{- 1,0,1\}^n:\ \e_j\neq
-1\}}\int_{\Z_m^n}d_\M(f(x+\e),f(x))^pd\mu(x)d\sigma(\e)\\
& \le
2^p\int_{\{-1,0,1\}^n}\int_{\Z_m^n}d_\M(f(x+\e),f(x))^pd\mu(x)d\sigma(\e).
\end{align*}
Summing over $j=1,\ldots,n$ yields the required result.
\end{proof}

\begin{lemma}\label{lem:with zeros} Let $(\M,d_\M)$ be a metric
space. Assume that for an integer $n$ and an even integer $m$, we have that for every $\ell\le n$, and every
$f:\Z_m^\ell \to \M$,
\begin{multline*}
\sum_{j=1}^\ell\int_{\Z_m^\ell}
d_\M\left(f\left(x+\frac{m}{2}e_j\right),f\left(x\right)\right)^pd\mu(x)\\ \le
C^pm^pn^{1-\frac{p}{q}}\Bigg(\E_\e\int_{\Z_m^\ell}d_\M\left(f(x+\e),f(x)\right)^pd\mu(x)\\
+\frac{1}{\ell}\sum_{j=1}^\ell\int_{\Z_m^\ell}
d_\M\left(f(x+e_j),f(x)\right)^pd\mu(x)\Bigg).
\end{multline*}
Then
$$
\Gamma_q^{(p)}(\M;n,m)\le 5C.
$$
\end{lemma}

\begin{proof}  Fix $f:\Z_m^n\to \M$ and $\emptyset \neq A\subseteq
\{1,\ldots,n\}$. Our assumption implies that
\begin{multline*}
 \sum_{j\in A}\int_{\Z_m^n}
d_\M\left(f\left(x+\frac{m}{2}e_j\right),f\left(x\right)\right)^pd\mu(x)\\[2pt]
  \le
C^pm^pn^{1-\frac{p}{q}}\Biggl(\E_\e\int_{\Z_m^n}d_\M\Biggl(f\Biggl(x+\sum_{j\in
A}\e_je_j\Biggr),f(x)\Biggr)^pd\mu(x)\\[2pt]
  +\frac{1}{|A|}\sum_{j\in
A}\int_{\Z_m^n} d_\M\left(f(x+e_j),f(x)\right)^pd\mu(x)\Biggr).
\end{multline*}
Multiplying this inequality by $\frac{2^{|A|}}{3^n}$, and summing
over all $\emptyset \ne A\subseteq \{1,\ldots,n\}$, we see that
\begin{eqnarray}&&\label{eq:manor's catch}\\[6pt]
&& \frac23 \sum_{j=1}^n\int_{\Z_m^n}
d_\M\left(f\left(x+\frac{m}{2}e_j\right),f\left(x\right)\right)^pd\mu(x)\nonumber\\[6pt]
&&\quad =\sum_{\emptyset\ne
A\subseteq \{1,\ldots,n\}}\frac{2^{|A|}}{3^n}\sum_{j\in
A}\int_{\Z_m^n}
d_\M\left(f\left(x+\frac{m}{2}e_j\right),f\left(x\right)\right)^pd\mu(x)
\nonumber\\[6pt]
&&\quad \le C^pm^pn^{1-\frac{p}{q}}\Biggl(\sum_{\emptyset
\ne A\subseteq
\{1,\ldots,n\}}\frac{2^{|A|}}{3^n}\E_\e\int_{\Z_m^n}d_\M\Biggl(f\Biggl(x+\sum_{j\in
A}\e_je_j\Biggr),f(x)\Biggr)^pd\mu(x)\nonumber \end{eqnarray}
\begin{eqnarray}&&\qquad +
\sum_{\emptyset \ne A\subseteq
\{1,\ldots,n\}}\frac{2^{|A|}}{|A|3^n}\sum_{j\in A}\int_{\Z_m^n}
d_\M\left(f(x+e_j),f(x)\right)^pd\mu(x)\Biggr)\nonumber\\[6pt] \label{eq:use
triangle}&&\quad \le
C^pm^pn^{1-\frac{p}{q}}\Bigg(\int_{\{-1,0,1\}^n}\int_{\Z_m^n}
d_\M\left(f\left(x+\delta\right),f(x)\right)^pd\mu(x)d\sigma(\delta)\\[6pt]
&&\nonumber\quad\phantom{\le
C^pm^pn^{1-\frac{p}{q}}\Bigg(}
+\frac{1}{n}\sum_{j=1}^n\int_{\Z_m^n}
d_\M\left(f(x+e_j),f(x)\right)^pd\mu(x)\Bigg)\\[6pt] &&\le
C^pm^pn^{1-\frac{p}{q}}\left(3^p+1\right)
\int_{\{-1,0,1\}^n}\int_{\Z_m^n}d_\M\left(f\left(x+\delta\right),f(x)\right)^pd\mu(x)d\sigma(\delta),\nonumber
\end{eqnarray}
where we used the fact that in~\eqref{eq:manor's catch}, the
coefficient of $d_\M\left(f(x+e_j),f(x)\right)^p$ equals
$\sum_{k=1}^{n}\frac{2^k}{k3^n}\binom{n-1}{k-1}\le \frac{1}{n}$,
and in~\eqref{eq:use triangle} we used Lemma~\ref{lem:pass to
diagonals}.
\end{proof}

\section{Warmup: the case of Hilbert space}

The fact that Hilbert spaces have metric cotype $2$ is
particularly simple to prove. This is contained in the following
proposition.

\begin{prop}\label{prop:hilbert}
Let $H$ be a Hilbert space. Then for every integer $n${\rm ,} and every
integer $m\ge \frac23\pi\sqrt{n}$ which is divisible by $4${\rm ,}
$$
\Gamma_2(H;n,m)\le \frac{\sqrt{6}}{\pi}.
$$
\end{prop}

\begin{proof}  Fix $f:Z_m^n\to H$ and decompose it into Fourier coefficients:
$$
f(x)=\sum_{k\in \mathbb Z_m^n} W_k(x)\widehat f(k).
$$
For every $j=1,2,\ldots,n$ we have that
$$
f\left(x+\frac{m}{2}e_j\right)-f(x)=\sum_{k\in \mathbb Z_m^n}
W_k(x)\left(e^{\pi ik_j}-1\right)\widehat f(k).
$$
Thus
\begin{eqnarray*}
&&\sum_{j=1}^n\int_{\Z_m^n}
\left\|f\left(x+\frac{m}{2}e_j\right)-f(x)\right\|_H^2d\mu(x) \\&=& \sum_{k\in
\Z_m^n}\Biggl(\sum_{j=1}^n\left|e^{\pi ik_j}-1
\right|^2\Biggr)\left\|\widehat f(k)\right\|_H^2 =  4\sum_{k\in
\Z_m^n}|\{j:   k_j\equiv 1\hbox{ mod}\, 2\}|\cdot\left\|\widehat
f(k)\right\|_H^2.
\end{eqnarray*}
Additionally, for every $\e\in \{-1,0,1\}^n$,
$$
f(x+\e)-f(x)=\sum_{k\in \mathbb Z_m^n} W_k(x)(W_k(\e)-1)\widehat
f(k).
$$
Thus
\begin{multline*}
\int_{\{-1,0,1\}^n}\int_{\Z_m^n}\|f(x+\e)-f(x)\|_H^2d\mu(x)d\sigma(\e)\\
=\sum_{k\in
\mathbb
Z_m^n}\left(\int_{\{-1,0,1\}^n}\left|W_k(\e)-1\right|^2d\sigma(\e)\right)
\left\|\widehat f(k)\right\|^2_H.
\end{multline*}
Observe that
\begin{eqnarray*}
\int_{\{-1,0,1\}^n}\left|W_k(\e)-1\right|^2d\sigma(\e)&=&\int_{\{-1,0,1\}^n}\Biggl|\exp\Biggl(\frac{2\pi
i}{m}\sum_{j=1}^m
k_j\e_j\Biggl)-1\Biggr|^2d\sigma(\e)\\
&=&2-2\, \mathrm{Re}\prod_{j=1}^n\int_{\{-1,0,1\}^n}\exp\left(\frac{2\pi
i }{m}k_j\e_j\right)d\sigma(\e)\\
&=& 2-2\prod_{j=1}^n\frac{1+2\cos\left(\frac{2\pi
}{m}k_j\right)}{3}\\
&\ge&
 2-2\prod_{j:\ k_j\equiv 1\mod 2}\frac{1+2\left|\cos\left(\frac{2\pi
}{m}k_j\right)\right|}{3}.
\end{eqnarray*}
Note that if $m$ is divisible by $4$ and $\ell\in
\{0,\ldots,m-1\}$ is an odd integer, then
$$
\left|\cos\left(\frac{2\pi\ell}{m}\right)\right|\le
\left|\cos\left(\frac{2\pi}{m}\right)\right|\le
1-\frac{\pi^2}{m^2}.
$$
Hence
\begin{eqnarray*}
\int_{\{-1,0,1\}^n} \left|W_k(\e)-1\right|^2 d \sigma(\e) &\ge&
2\Biggl(1-\Biggl(1-\frac{2\pi^2}{3m^2}\Biggr)^{|\{j:\ k_j\equiv
1\mod 2\}|}\Biggr)\\&\ge& 2\Biggl(1-e^{-\frac{2|\{j:\ k_j\equiv
1\mod 2\}|\pi^2}{3m^2}}\Biggr)\\&\ge& |\{j:\ k_j\equiv 1\mod
2\}|\cdot \frac{2\pi^2}{3m^2},
\end{eqnarray*}
provided that $m\ge \frac23\pi \sqrt{n}$.
\end{proof}

\section{$K$-convex spaces}

In this section we prove the ``hard direction" of
Theorem~\ref{thm:cotype} and Theorem~\ref{thm:weak cotype} when
$X$ is a $K$-convex Banach space; namely, we show that in this case
Rademacher cotype $q$ implies metric cotype $q$. There are two
reasons why we single out this case before passing to the proofs
of these theorems in full generality. First of all, the proof for
$K$-convex spaces is different and simpler than the general case.
More importantly, in the case of $K$-convex spaces we are able to
obtain optimal bounds on the value of $m$ in
Definition~\ref{def:cotype} and Definition~1.3.
Namely, we show that if $X$ is a $K$-convex Banach space of cotype
$q$, then for every $1\le p\le q$,
$m_q^{(p)}(X;n,\Gamma)=O(n^{1/q})$, for some $\Gamma=\Gamma(X)$.
This is best possible due to Lemma~\ref{lem:lower m}. In the case
of general Banach spaces we obtain worse bounds, and this is why
we have the restriction that $X$ is $K$-convex in
Theorem~\ref{thm:uniform} and Theorem~\ref{thm:coarse}. This issue
is taken up again in Section~\ref{section:problems}.

\begin{theorem}\label{thm:K} Let $X$ be a $K$-convex Banach space with cotype $q$. Then for
every integer $n$ and every integer $m$ which is divisible by $4${\rm ,}
$$
m\ge \frac{2n^{1/q}}{C_q^{(p)}(X)K_p(X)} \implies
\Gamma_q^{(p)}(X;n,m)\le 15C_q^{(p)}(X)K_p(X).
$$
\end{theorem}

\begin{proof} For $f:\Z_m^n\to X$ we define the following operators:
\begin{eqnarray*}
\widetilde \partial_j f(x)&=&f(x+e_j)-f(x-e_j),
\\
\mathcal E_j f(x)&=&\E_\e f\Biggl(x+\sum_{\ell\neq j} \e_\ell
e_\ell\Biggr),
\end{eqnarray*}
and for $\e\in \{-1,0,1\}^n$,
$$
\partial_\e f(x)=f(x+\e)-f(x).  
$$
These operators operate diagonally on the Walsh basis
$\{W_k\}_{k\in \Z_m^n}$ as follows:
\begin{eqnarray}\label{eq:partial tilde}
\widetilde \partial_j W_k&= &\left(W_k(e_j)-W_k(-e_j)\right)W_k=
2\sin\left(\frac{2\pi i k_j}{m}\right)\cdot W_k,
\\
\label{eq:Ej}
\mathcal E_j W_k &=&\Biggl(\E_\e \prod_{\ell\neq j}e^{\frac{2\pi i
\e_\ell k_\ell}{m}}\Biggr)W_k=\Biggl(\prod_{\ell\neq j}
\cos\left(\frac{2\pi k_\ell}{m}\right)\Biggr)W_k,
\end{eqnarray}
and for $\e\in \{-1,1\}^n$,
\begin{eqnarray}\label{eq:partial}
\partial_\e
W_k&=&\left(W(\e)-1\right)W_k\\\nonumber&=&\Biggl(\prod_{j=1}^n
e^{\frac{2\pi i \e_j k_j}{m}}-1\Biggr)W_k\\&=&\nonumber
\Biggl(\prod_{j=1}^n\Biggl(\cos\left(\frac{2\pi 
\e_jk_j}{m}\right)+i\sin\left(\frac{2\pi 
\e_jk_j}{m}\right)\Biggr)-1\Biggr)W_k\\&=&
\Biggl(\prod_{j=1}^n\left(\cos\left(\frac{2\pi 
k_j}{m}\right)+i\e_j\sin\left(\frac{2\pi 
k_j}{m}\right)\right)-1\Biggr)W_k.\nonumber
\end{eqnarray}
The last step was a crucial observation, using the fact that
$\e_j\in \{-1,1\}$. Thinking of $\partial_\e W_k$ as a function of
$\e\in \{-1,1\}^n$, equations~\eqref{eq:partial tilde},
\eqref{eq:Ej} and~\eqref{eq:partial} imply that
\begin{eqnarray*}
\Rad(\partial_\e W_k)&=&i\Biggl(\sum_{j=1}^n
\e_j\sin\left(\frac{2\pi k_j}{m}\right)\cdot \prod_{\ell\neq j}
\cos\left(\frac{2\pi
k_\ell}{m}\right)\Biggr)W_k\\
&=&\frac{i}{2}\Biggl(\sum_{j=1}^n \e_j
\widetilde \partial_j \mathcal E_j\Biggr)W_k.
\end{eqnarray*}
Thus for every $x\in \Z_m^n$ and $f:\Z_m^n\to X$,
$$
\Rad(\partial_\e f(x))=\frac{i}{2}\Biggl(\sum_{j=1}^n \e_j
\widetilde
\partial_j \mathcal E_j\Biggr)f(x).
$$
It follows that
\begin{eqnarray}\label{eq:rademacher case}
&&\hskip-16pt \int_{\Z_m^n} \E_\e\Biggl\|\sum_{j=1}^n \e_j
\bigl[ \mathcal E_j
f(x+e_j)-\mathcal E_j f(x-e_j) \bigr ]\Biggr\|_X^pd\mu(x)\\
 &&\qquad\qquad = \int_{\Z_m^n}
\E_\e\Biggl\|\sum_{j=1}^n\e_j \widetilde \partial_j \mathcal E_j
f(x) \Biggr\|_X^pd\mu(x)\nonumber\\&&\qquad\qquad =\int_{\Z_m^n}
\E_\e\|\Rad(\partial_\e f(x))\|_X^pd\mu(x)\nonumber\\ &&\qquad\qquad  \le 
K_p(X)^p\int_{\Z_m^n} \E_\e\|\partial_\e f(x)\|_X^pd\mu(x).
 \nonumber
\end{eqnarray}
By~\eqref{eq:rademacher case} and the definition of
$C_q^{(p)}(X)$, for every $C>C^{(p)}_q(X)$ we have that
\begin{eqnarray}\label{eq:holder}
&&
[K_p(X)C]^p\E_\e\int_{\Z_m^n}
\|f(x+\e)-f(x)\|_X^pd\mu(x)\\\nonumber&&\qquad \ge  C^p\cdot\E_\e
\int_{\Z_m^n}\Biggl\|\sum_{j=1}^n \e_j[\mathcal E_j
f(x+e_j)-\mathcal E_j
f(x-e_j)]\Biggr\|_X^pd\mu(x)\\&&\qquad \ge\int_{\Z_m^n}\Biggl(\sum_{j=1}^n
\left\|\mathcal E_j f(x+e_j)-\mathcal E_j
f(x-e_j)\right\|_X^q\Biggr)^{p/q}d\mu(x)\nonumber\\&&\qquad
\ge\frac{1}{n^{1-p/q}}\sum_{j=1}^n\int_{\Z_m^n}\left\|\mathcal E_j f(x+e_j)-\mathcal E_j
f(x-e_j)\right\|_X^pd\mu(x). \nonumber
\end{eqnarray}
Now, for $j\in \{1,\ldots,n\}$,
\begin{eqnarray}\label{eq:after integral}
\qquad &&\int_{\Z_m^n} \left\|\mathcal
E_jf\left(x+\frac{m}{2}e_j\right)-\mathcal
E_jf\left(x\right)\right\|_X^pd\mu(x)\\
&&\qquad \le\left(\frac{m}{4}\right)^{p-1}
\sum_{s=1}^{m/4}\int_{\Z_m^n} \left\|\mathcal
E_jf\left(x+2se_j\right)-\mathcal
E_jf\left(x+2(s-1)e_j\right)\right\|_X^pd\mu(x)\nonumber\\&&\qquad =
\left(\frac{m}{4}\right)^p\int_{\Z_m^n} \left\|\mathcal E_j
f(x+e_j)-\mathcal E_j f(x-e_j)\right\|_X^pd\mu(x).\nonumber
\end{eqnarray}
Plugging~\eqref{eq:after integral} into~\eqref{eq:holder} we get
\begin{eqnarray*}
&&  \hskip-38pt
\left(\frac{m}{4}\right)^{p}n^{1-\frac{p}{q}}[K_p(X)C]^p\E_\e\int_{\Z_m^n}
\|f(x+\e)-f(x)\|_X^pd\mu(x)\\
&&\quad \ge \sum_{j=1}^n\int_{\Z_m^n}
\left\|\mathcal E_jf\left(x+\frac{m}{2}e_j\right)-\mathcal
E_jf\left(x\right)\right\|_X^pd\mu(x)\\ &&\quad \ge\frac{1}{3^{p-1}}\sum_{j=1}^n\int_{\Z_m^n}
\left\|f\left(x+\frac{m}{2}e_j\right)-f\left(x\right)\right\|_X^pd\mu(x)\\
&&\qquad -2\sum_{j=1}^n\int_{\Z_m^n}
\left\|\mathcal
E_jf\left(x\right)-f\left(x\right)\right\|_X^pd\mu(x)\\&&\quad
=\frac{1}{3^{p-1}}\sum_{j=1}^n\int_{\Z_m^n}
\left\|f\left(x+\frac{m}{2}e_j\right)-f\left(x\right)\right\|_X^pd\mu(x)\\
&&\qquad -2\sum_{j=1}^n\int_{\Z_m^n}
\Biggl\|\E_\e\Biggl(f\Biggl(x+\sum_{\ell\neq j} \e_\ell
e_\ell\Biggr)-f\left(x\right)\Biggr)\Biggr\|_X^pd\mu(x)\nonumber\\&&\quad \ge
\frac{1}{3^{p-1}}\sum_{j=1}^n\int_{\Z_m^n}
\left\|f\left(x+\frac{m}{2}
e_j\right)-f\left(x\right)\right\|_X^pd\mu(x)\\*
&&\qquad -2\sum_{j=1}^n\E_\e\int_{\Z_m^n}
\Biggl\|f\Biggl(x+\sum_{\ell\neq j} \e_\ell
e_\ell\Biggr)-f\left(x\right)\Biggr\|_X^pd\mu(x)\\
&&\quad \ge \frac{1}{3^{p-1}}\sum_{j=1}^n\int_{\Z_m^n}
\left\|f\left(x+\frac{m}{2}
e_j\right)-f\left(x\right)\right\|_X^pd\mu(x)\\
&&\qquad -2^pn\E_\e\int_{\Z_m^n}
\left\|f\left(x+\e\right)-f\left(x\right)\right\|_X^pd\mu(x)\\& &\qquad -
2^p\sum_{j=1}^n\E_\e\int_{\Z_m^n}
\left\|f\left(x+\e_je_j\right)-f\left(x\right)\right\|_X^pd\mu(x).
\end{eqnarray*}
Thus, the required result follows from Lemma~\ref{lem:with zeros}.
\end{proof}

 The above argument actually gives the following
generalization of Theorem~\ref{thm:K}, which holds for products of
arbitrary compact Abelian groups.

\begin{theorem}\label{thm:groups}Let $G_1,\ldots,G_n$ be compact
Abelian groups{\rm ,} $(g_1,\ldots,g_n)\in G_1\times\cdots \times G_n${\rm ,}
and let $X$ be a $K$-convex Banach space. Then for every integer
$k$ and every $f:G_1\times\cdots\times G_n\to X${\rm ,}
\vglue-18pt

\begin{small}
\begin{multline*}
\sum_{j=1}^n \int_{G_1\times \cdots\times G_n}
\|f(x+2kg_je_j)-f(x)\|_X^pd(\mu_{G_1}\otimes\cdots
\otimes\mu_{G_n})(x)\\\le C^p \int_{\{-1,0,1\}^n}\int_{G_1\times
\cdots\times
G_n}\Le\|f\Le(x+\sum_{j=1}^n\e_jg_je_j\Ri)-f(x)\Ri\|_X^pd(\mu_{G_1}\otimes\cdots
\otimes\mu_{G_n})(x)d\sigma(\e),
\end{multline*} \end{small}

\vglue-8pt\noindent 
where
$$
C\le
5\max\left\{C_q^{(p)}(X)K_p(X)kn^{\frac{1}{p}-\frac{1}{q}},n^{\frac{1}{p}}\right\}.
$$
\end{theorem}

Here $\mu_G$ denotes the normalized Haar measure on a compact
Abelian group $G$. We refer the interested reader to the
book~\cite{Rudin90}, which contains the necessary background
required to generalize the proof of Theorem~\ref{thm:K} to this
setting.

\section{The equivalence of Rademacher cotype and metric cotype}

We start by establishing the easy direction in
Theorem~\ref{thm:cotype} and Theorem~\ref{thm:weak cotype}, i.e.
that metric cotype implies Rademacher cotype.

\subsection{Metric cotype implies Rademacher
cotype}\label{section:easy direction}
Let $X$ be a Banach space and assume that
$\Gamma_q^{(p)}(X)<\infty$ for some $1\le p\le q$. Fix
$\Gamma>\Gamma_q^{(p)}(X)$, $v_1,\ldots,v_n\in X$, and let $m$ be
an even integer. Define $f:\Z_m^n\to X$ by
$$
f(x_1,\ldots,x_n)=\sum_{j=1}^n e^{\frac{2\pi i x_j}{m}}v_j.
$$
Then
\begin{eqnarray}\label{eq:reverse cotype}
\sum_{j=1}^n\int_{\Z_m^n}\left\|f\left(x+\frac{m}{2}e_j\right)-f(x)\right\|_X^pd\mu(x)=2^p\sum_{j=1}^n
\|v_j\|_X^p,
\end{eqnarray}
and
\begin{multline}\label{eq:reverse cotype2}
\int_{\{-1,0,1\}^n}\int_{\Z_m^n}\left\|f\left(x+\delta\right)-f(x)\right\|_X^pd\mu(x)d\sigma(\delta)\\
=
\int_{\{-1,0,1\}^n}\int_{\Z_m^n}\Biggl\|\sum_{j=1}^n e^{\frac{2\pi
i x_j}{m}}\left(e^{\frac{2\pi i \delta_j}{m}}-1\right)v_j
\Biggr\|_X^pd\mu(x)d\sigma(\delta).
\end{multline}

We recall the {\em contraction principle} (see~\cite{LT91}), which
states that for every $a_1,\ldots, a_n\in \R$,
$$
\E_\e\Biggl\|\sum_{j=1}^n \e_ja_jv_j\Biggr\|_X^p\le
\left(\max_{1\le j\le n} |a_j|\right)^p\cdot
\E_\e\Biggl\|\sum_{j=1}^n \e_jv_j\Biggr\|_X^p.
$$

Observe that for every $\e=(\e_1,\ldots,\e_n)\in \{-1,1\}^n$,

\begin{eqnarray*}
&&\!\!\!\!\!\!\!\!\!\!\!\!\!\!\!\!\!\!\!\!\int_{\{-1,0,1\}^n}\int_{\Z_m^n}\Le\|\sum_{j=1}^n
e^{\frac{2\pi i x_j}{m}}\left(e^{\frac{2\pi i
\delta_j}{m}}-1\right)v_j
\Ri\|_X^pd\mu(x)d\sigma(\delta)\\&=&\int_{\{-1,0,1\}^n}\int_{\Z_m^n}\Le\|\sum_{j=1}^n
e^{\frac{2\pi i}{m}
\left(x_j+\frac{m(1-\e_j)}{4}\right)}\left(e^{\frac{2\pi i
\delta_j}{m}}-1\right)v_j
\Ri\|_X^pd\mu(x)d\sigma(\delta)\\
&=&\int_{\{-1,0,1\}^n}\int_{\Z_m^n}\Le\|\sum_{j=1}^n \e_j
e^{\frac{2\pi i x_j}{m}}\left(e^{\frac{2\pi i
\delta_j}{m}}-1\right)v_j \Ri\|_X^pd\mu(x)d\sigma(\delta).
\end{eqnarray*}

Taking expectation with respect to $\e$, and using the contraction
principle, we see that
\begin{eqnarray}\label{eq:reverse cotype 2nd}
&& \int_{\{-1,0,1\}^n}\int_{\Z_m^n}\Le\|\sum_{j=1}^n
e^{\frac{2\pi i x_j}{m}}\left(e^{\frac{2\pi i
\delta_j}{m}}-1\right)v_j
\Ri\|_X^pd\mu(x)d\sigma(\delta)\\&&\quad=
\int_{\{-1,0,1\}^n}\int_{\Z_m^n}\E_\e\Le\|\sum_{j=1}^n \e_j
e^{\frac{2\pi i x_j}{m}}\left(e^{\frac{2\pi i
\delta_j}{m}}-1\right)v_j
\Ri\|_X^pd\mu(x)d\sigma(\delta)\nonumber\\&&\quad \le
\int_{\{-1,0,1\}^n}\int_{\Z_m^n}2^p\left(\max_{1\le j\le
n}\left|e^{\frac{2\pi i
\delta_j}{m}}-1\right|\right)^p\E_\e\Le\|\sum_{j=1}^n \e_j v_j
\Ri\|_X^pd\mu(x)d\sigma(\delta)\nonumber\\*&&\quad
\le\left(\frac{4\pi}{m}\right)^p\E_\e\Le\|\sum_{j=1}^n \e_j v_j
\Ri\|_X^p,\nonumber
\end{eqnarray}
where in the last inequality above we used the fact that for
$\theta\in [0,\pi]$, $|e^{i\theta}-1|\break\le \theta$.

Combining~\eqref{eq:def weak cotype}, \eqref{eq:reverse cotype},
\eqref{eq:reverse cotype2}, and~\eqref{eq:reverse cotype 2nd}, we
get that
\begin{eqnarray*}
2^p\sum_{j=1}^n \|v_j\|_X^p\le
\Gamma^pm^p\left(\frac{4\pi}{m}\right)^pn^{1-\frac{p}{q}}\E_\e\Le\|\sum_{j=1}^n
\e_jv_j
\Ri\|_X^p\!\!=\!\left(4\pi\Gamma\right)^pn^{1-\frac{p}{q}}\E_\e\Le\|\sum_{j=1}^n
\e_jv_j \Ri\|_X^p.
\end{eqnarray*}
If $p=q$ we see that $C_q(X)\le 2\pi \Gamma_q(X)$. If $p<q$ then
when $\|v_1\|_X=\cdots=\|v_n\|_X=1$ we get that
$$
\Le(\E_\e\Le\|\sum_{j=1}^n \e_jv_j \Ri\|_X^q\Ri)^{1/q}\ge
\Le(\E_\e\Le\|\sum_{j=1}^n \e_jv_j
\Ri\|_X^p\Ri)^{1/p}=\Omega\left(\frac{n^{1/q}}{\Gamma}\right).
$$
This means that $X$ has ``equal norm cotype $q$", implying that
$X$ has cotype $q'$ for every $q'>q$ (see~\cite{Tza79}, \cite{KT81}, \cite{T-J89}
for quantitative versions of this statement). When $q=2$ this
implies that $X$ has cotype $2$ (see~\cite{Tza79} and the
references therein).
 
\subsection{Proof of Theorem~{\rm \ref{thm:cotype}} and Theorem~{\rm \ref{thm:weak cotype}}}
The proof of Theorem~\ref{thm:cotype} and Theorem~\ref{thm:weak
cotype} is based on several lemmas. Fix an odd integer $k\in
\mathbb N$, with $k< \frac{m}{2}$, and assume that $1\le p\le q$.
Given $j\in \{1,\ldots,n\}$, define $S(j,k)\subseteq \Z_m^n$ by
$$
S(j,k):= \left\{y\in [-k,k]^n\subseteq \Z_m^n:\ y_j\equiv
0\mod 2\ \mathrm{and}\ \forall\  \ell\neq j,\ y_\ell\equiv 1\mod
2\right\}.
$$
For $f:\Z_m^n\to X$ we define
\begin{eqnarray}\label{eq:def Ek}
\op_j^{(k)}f(x)=\left(f*\frac{\mathbf{1}_{S(j,k)}}{\mu(S(j,k))}\right)(x)=\frac{1}{\mu(S(j,k))}\int_{S(j,k)}
f(x+y)d\mu(y).
\end{eqnarray}

\begin{lemma}\label{lem:approx} \hskip-5pt
For every $p\ge 1${\rm ,} every $j\in \{1,\ldots,n\}${\rm ,} and every
$f:\Z_m^n\to X${\rm ,}
\begin{eqnarray*}
\int_{\Z_m^n}\left\|\op^{(k)}_jf(x)-f(x)\right\|_X^pd\mu(x)& \le&
2^pk^p
\E_\e\int_{\Z_m^n}\|f(x+\e)-f(x)\|_X^pd\mu(x)\\&&+2^{p-1}\int_{\Z_m^n}
\|f(x+e_j)-f(x)\|_X^pd\mu(x).
\end{eqnarray*}
\end{lemma}

\begin{proof}
By convexity, for every $x\in \Z_m^n$,
\begin{eqnarray}\label{eq:convexity}\qquad
\left\|\op_j^{(k)}f(x)-f(x)\right\|_X^p&=&\left\|\frac{1}{\mu(S(j,k))}\int_{S(j,k)}
[f(x+y)-f(x)]d\mu(y)\right\|_X^p\\&\le&
\frac{1}{\mu(S(j,k))}\int_{S(j,k)} \|f(x)-f(x+y)\|_X^pd\mu(y).\nonumber
\end{eqnarray}

Let $x\in \{0,\ldots,k\}^n$ be such that for all $j\in \{1,\ldots,n\}$,
$x_j$ is a positive odd integer. Observe that there exists a
geodesic $\gamma:\{0,1,\ldots,\|x\|_\infty\}\to \Z_m^n$ such that
$\gamma(0)=0$, $\gamma(\|x\|_\infty)=x$ and for every
$t\in\{1,\ldots,\|x\|_\infty\}$,\break $\gamma(t)-\gamma(t-1)\in
\{-1,1\}^n$. Indeed, we define $\gamma(t)$ inductively as follows:
$\gamma(0)=0$, $\gamma(1)=(1,1,\ldots,1)$, and if $t\ge 2$ is odd
then
$$
\gamma(t)=\gamma(t-1)+\sum_{s=1}^n e_s\quad \mathrm{and}\quad
\gamma(t+1)= \gamma(t-1)+2\sum_{\substack{s\in \{1,\ldots, n\}\\
\gamma(t-1)_s<x_s}} e_s.
$$
Since all the coordinates of $x$ are odd,
$\gamma(\|x\|_\infty)=x$. In what follows we fix an arbitrary
geodesic $\gamma_x:\{0,1,\ldots,\|x\|_\infty\}\to \Z_m^n$ as
above. For $x\in (\Z_m\setminus \{0\})^n$ we denote
$|x|=(|x_1|,\ldots,|x_n|)$ and
$\sign(x)=(\sign(x_1),\ldots,\sign(x_n))$. If $x\in [-k,k]^n$ is
such that all of its coordinates are odd, then we define
$\gamma_x=\sign(x)\cdot \gamma_{|x|}$ (where the multiplication is
coordinate-wise).

If $y\in S(j,k)$ then all the coordinates of $y\pm e_j$ are odd.
We can thus define two geodesic paths
$$
\gamma_{x,y}^{+1}=x+e_j+\gamma_{y-e_j}\quad \mathrm{and}\quad
\gamma_{x,y}^{-1}=x-e_j+\gamma_{y+e_j},
$$
where the addition is point-wise.

For $z\in \Z_m^n$ and $\e\in \{-1,1\}^n$ define
\begin{multline*}
F^{+1}(z,\e)=\Big\{(x,y)\in \Z_m^n\times S(j,k):\ \exists t\in
\{1,\ldots,\|y-e_j\|_\infty\},\\ \gamma^{+1}_{x,y}(t-1)=z,\
\gamma^{+1}_{x,y}(t)=z+\e\Big\},
\end{multline*}
and
\begin{multline*}
F^{-1}(z,\e)=\Big\{(x,y)\in \Z_m^n\times S(j,k):\ \exists t\in
\{1,\ldots,\|y+e_j\|_\infty\},\\ \gamma^{-1}_{x,y}(t-1)=z,\
\gamma^{-1}_{x,y}(t)=z+\e\Big\}.
\end{multline*}

\begin{claim}\label{claim:independence} For every $z,w\in \Z_m^n$ and $\e,\delta\in
\{-1,1\}^n${\rm ,}
$$
|F^{+1}(z,\e)|+|F^{-1}(z,\e)|=|F^{+1}(w,\delta)|+|F^{-1}(w,\delta)|.
$$
\end{claim}
\vglue8pt

\begin{proof} 
Define $\psi:\Z_m^n\times S(j,k)\to \Z_m^n\times S(j,k)$ by
$$
\psi(x,y)=(w-\e\delta z+\e\delta x,\e\delta y).
$$
We claim that $\psi$ is a bijection between $F^{+1}(z,\e)$ and
$F^{\e_j\delta_j}(w,\delta)$, and also $\psi$ is a bijection
between $F^{-1}(z,\e)$ and $F^{-\e_j\delta_j}(w,\delta)$. Indeed,
if $(x,y)\in F^{+1}(z,\e)$ then there exists $t\in
\{1,\ldots,\|y-e_j\|_\infty\}$ such that
$\gamma^{+1}_{x,y}(t-1)=z$ and $\gamma^{+1}_{x,y}(t)=z+\e$. The
path $w-\e\delta z+\e\delta \gamma_{x,y}^{+1}$ equals the path
$\gamma^{\e_j\delta_j}_{\psi(x,y)}$, which by definition goes
through $w$ at time $t-1$ and $w+\delta$ at time $t$. Since these
transformations are clearly invertible, we obtain the required
result for $F^{+1}(z,\e)$. The proof for $F^{-1}(z,\e)$ is
analogous.
\end{proof}

\begin{claim}\label{claim:fubini}
Denote $N=|F^{+1}(z,\e)|+|F^{-1}(z,\e)|${\rm ,} which is independent of
$z\in \Z_m^n$ and $\e\in \{-1,1\}^n${\rm ,} by
Claim~{\rm \ref{claim:independence}}. Then
$$
N\le \frac{k\cdot |S(j,k)|}{2^{n-1}}.
$$
\end{claim}

\begin{proof}
We have that
\begin{eqnarray*}
N\cdot m^n\cdot 2^n&=&\sum_{(z,\e)\in \Z_m^n\times
\{-1,1\}^n}\left(|F^{+1}(z,\e)|+|F^{-1}(z,\e)|\right)\\[6pt]&=&\hskip-3pt\sum_{(z,\e)\in
\Z_m^n\times \{-1,1\}^n}\hskip-3pt\Le(\sum_{(x,y)\in \Z_m^n\times
S(j,k)}\sum_{t=1}^{\|y-e_j\|_\infty}\mathbf
{1}_{\{\gamma_{x,y}^{+1}(t-1)=z\ \wedge\
\gamma_{x,y}^{+1}(t)=z+\e\}}\hskip-3pt\Ri)\\[6pt]& &+\hskip-3pt\sum_{(z,\e)\in
\Z_m^n\times \{-1,1\}^n}\hskip-3pt\Le(\sum_{(x,y)\in \Z_m^n\times
S(j,k)}\sum_{t=1}^{\|y+e_j\|_\infty}\mathbf
{1}_{\{\gamma_{x,y}^{-1}(t-1)=z\ \wedge\
\gamma_{x,y}^{-1}(t)=z+\e\}}\hskip-3pt\Ri)\\[6pt] &=&\sum_{(x,y)\in \Z_m^n\times
S(j,k)}\|y-e_j\|_\infty+ \sum_{(x,y)\in \Z_m^n\times
S(j,k)}\|y+e_j\|_\infty\\[6pt]
&\le& 2k \cdot m^n\cdot |S(j,k)|.
\end{eqnarray*}
\end{proof}

We now conclude the proof of Lemma~\ref{lem:approx}. Observe that
for $x\in \Z_m^n$ and $y\in S(j,k)$,
\begin{eqnarray}\label{eq:geodesic+}&&\\[-2pt] \nonumber
\frac{\|f(x)-f(x+y)\|_X^p}{2^{p-1}}&\le&
\|f(x)-f(x+e_j)\|_X^p\\
&&\nonumber+
\|y-e_j\|_\infty^{p-1}\sum_{t=1}^{\|y-e_j\|_\infty}
\|f(\gamma_{x,y}^{+1}(t))-f(\gamma_{x,y}^{+1}(t-1))\|_X^p\\
&\le& \|f(x)-f(x+e_j)\|_X^p\nonumber\\
&&+
k^{p-1}\sum_{t=1}^{\|y-e_j\|_\infty}
\|f(\gamma_{x,y}^{+1}(t))-f(\gamma_{x,y}^{+1}(t-1))\|_X^p,\nonumber
\\[-9pt]
\noalign{\noindent and}
 \label{eq:geodesic-}&&\\[-2pt]
\nonumber\frac{\|f(x)-f(x+y)\|_X^p}{2^{p-1}}&\le&
\|f(x)-f(x-e_j)\|_X^p\\
&&+
\|y+e_j\|_\infty^{p-1}\sum_{t=1}^{\|y+e_j\|_\infty}
\|f(\gamma_{x,y}^{-1}(t))-f(\gamma_{x,y}^{-1}(t-1))\|_X^p\nonumber\\
&\le& \|f(x)-f(x-e_j)\|_X^p\nonumber\\ &&+
k^{p-1}\sum_{t=1}^{\|y+e_j\|_\infty}
\|f(\gamma_{x,y}^{-1}(t))-f(\gamma_{x,y}^{-1}(t-1))\|_X^p.\nonumber
\end{eqnarray}
Averaging inequalities~\eqref{eq:geodesic+}
and~\eqref{eq:geodesic-}, and integrating, we get that
\begin{eqnarray}\label{eq:use N}
&&\frac{1}{\mu(S(j,k))}\int_{\Z_m^n}\int_{S(j,k)}
\|f(x)-f(x+y)\|_X^pd\mu(y)d\mu(x)\\&&\quad \le \nonumber
2^{p-1}\int_{\Z_m^n}\|f(x+e_j)-f(x)\|_X^pd\mu(x)\nonumber\\
&&\nonumber\qquad +
(2k)^{p-1}\frac{N\cdot
2^{n}}{|S(j,k)|}\E_\e\int_{\Z_m^n}\|f(z+\e)-f(z)\|_X^pd\mu(z)\\ \label{eq:use
bound on N} &&\quad \le
2^{p-1}\int_{\Z_m^n}\|f(x+e_j)-f(x)\|_X^pd\mu(x)\\
&&\qquad +
(2k)^{p}\E_\e\int_{\Z_m^n}\|f(z+\e)-f(z)\|_X^pd\mu(z),\nonumber
\end{eqnarray}
where in~\eqref{eq:use N} we used Claim~\ref{claim:independence}
and in~\eqref{eq:use bound on N} we used Claim~\ref{claim:fubini}.
By~\eqref{eq:convexity}, this completes the proof of
Lemma~\ref{lem:approx}.
\end{proof}

Lemma~\ref{lem:pass to average on e} below is the heart of our
proof. It contains the cancellation of terms which is key to the
validity of Theorem~\ref{thm:cotype} and Theorem~\ref{thm:weak
cotype}.

\begin{lemma}\label{lem:pass to average on e} For every $f:\Z_m^n\to X${\rm ,} 
every integer $n${\rm ,} every even integer~$m${\rm ,} every $\e\in \{-1,1\}^n${\rm ,}
 every odd integer $k<m/2$,
and every $p\ge 1${\rm ,}
\begin{eqnarray*}&&\hskip-30pt
\int_{\Z_m^n} \Le\|\sum_{j=1}^n
\e_j\left[\op_j^{(k)}f(x+e_j)-\op_j^{(k)}f(x-e_j)\right]\Ri\|_X^pd\mu(x)\\
&&\qquad \le
3^{p-1}
\int_{\Z_m^n}\|f(x+\e)-f(x-\e)\|_X^pd\mu(x)\\
&&\qquad\quad +\frac{\cdot24^pn^{2p-1}}{k^{p}}\sum_{j=1}^n\int_{\Z_m^n}
\|f(x+e_j)-f(x)\|_X^pd\mu(x) .
\end{eqnarray*}
\end{lemma}

We postpone the proof of Lemma~\ref{lem:pass to average on e} to
Section~\ref{section:the lemma}, and proceed to prove
Theorem~\ref{thm:cotype} and Theorem~\ref{thm:weak cotype}
assuming its validity.

\begin{proof}[Proof of Theorem~{\rm \ref{thm:cotype}} and Theorem~{\rm \ref{thm:weak cotype}}]
 Taking
expectations with respect to $\e\in \{-1,1\}^n$ in
Lemma~\ref{lem:pass to average on e} we get that
\begin{eqnarray}\label{eq:E}
&& \\[-4pt]
&&\E_\e\int_{\Z_m^n} \Le\|\sum_{j=1}^n
\e_j\left[\op_j^{(k)}f(x+e_j)-\op_j^{(k)}f(x-e_j)\right]\Ri\|_X^pd\mu(x)\nonumber\\
& &\quad\le 3^{p-1}\E_\e
\int_{\Z_m^n}2^{p-1}\left(\|f(x+\e)-f(x)\|_X^p+\|f(x)-f(x-\e)\|_X^p\right)d\mu(x)\nonumber\\
&&\qquad +
\frac{24^pn^{2p-1}}{k^{p}}\sum_{j=1}^n\int_{\Z_m^n}
\|f(x+e_j)-f(x)\|_X^pd\mu(x)\nonumber\\&&\quad \le
\frac{6^p}{3}\E_\e\int_{\Z_m^n}\|f(x+\e)-f(x)\|_X^pd\mu(x)\nonumber\\
&&\qquad +\frac{24^pn^{2p-1}}{k^{p}}\sum_{j=1}^n\int_{\Z_m^n}
\|f(x+e_j)-f(x)\|_X^pd\mu(x).\nonumber
\end{eqnarray}

Fix $x\in \Z_m^n$ and let $m$ be an integer which is divisible by
$4$ such that $m\ge 6n^{2+1/q}$. Fixing $C>C_q^{(p)}(X)$, and
applying  the definition of $C_q^{(p)}(X)$ to the vectors
$\left\{\op_j^{(k)} f(x+e_j)-\op_j^{(k)}
f(x-e_j)\right\}_{j=1}^n$, we get
\begin{multline}\label{eq:just used cotype}
\E_\e\Le\|\sum_{j=1}^n \e_j\left[\op_j^{(k)} f(x+e_j)-\op_j^{(k)}
f(x-e_j)\right]\Ri\|_X^p\\
\ge \frac{1}{C^p\cdot
n^{1-p/q}}\sum_{j=1}^n \left\|\op_j^{(k)} f(x+e_j)-\op_j^{(k)}
f(x-e_j)\right\|_X^p.
\end{multline}
Now, for every $j\in \{1,\ldots,n\}$,
\begin{multline}\label{eq:before average}
\sum_{s=1}^{m/4}
\left\|\op_j^{(k)}f\left(x+2se_j\right)-\op_j^{(k)}f\left(x+2(s-1)e_j\right)\right\|_X^p\\
\ge
\left(\frac{4}{m}\right)^{p-1}\left\|\op_j^{(k)}f\left(x+\frac{m}{2}e_j\right)-\op_j^{(k)}f\left(x\right)\right\|_X^p.
\end{multline}
Averaging~\eqref{eq:before average} over $x\in \Z_m^n$ we get that
\begin{multline}\label{eq:after average}
\int_{\Z_m^n}
\left\|\op_j^{(k)}f(x+e_j)-\op_j^{(k)}f(x-e_j)\right\|_X^pd\mu(x)\\
\ge\left(\frac{4}{m}\right)^p\int_{\Z_m^n}
\left\|\op_j^{(k)}f\left(x+\frac{m}{2}e_j\right)-\op_j^{(k)}f\left(x\right)\right\|_X^pd\mu(x).
\end{multline}
Combining~\eqref{eq:just used cotype} and~\eqref{eq:after average}
we get the inequality
\begin{multline}\label{eq:before lemma}
\E_\e \int_{\Z_m^n}\Le\|\sum_{j=1}^n \e_j\left[\op_j^{(k)}
f(x+e_j)-\op_j^{(k)} f(x-e_j)\right]\Ri\|_X^pd\mu(x)\\\ge
\frac{1}{C^p\cdot n^{1-p/q}}\cdot
\left(\frac{4}{m}\right)^p\sum_{j=1}^n\int_{\Z_m^n}
\left\|\op_j^{(k)}f\left(x+\frac{m}{2}e_j\right)-\op_j^{(k)}f\left(x\right)\right\|_X^pd\mu(x).
\end{multline}
Now, for every $j\in \{1,\ldots,n\}$,
\begin{eqnarray}\label{eq:each j}
\qquad&&\hskip-35pt \int_{\Z_m^n}
\left\|\op_j^{(k)}f\left(x+\frac{m}{2}e_j\right)-\op_j^{(k)}f\left(x\right)\right\|_X^pd\mu(x)\\
&&\quad \ge\frac{1}{3^{p-1}}\int_{\Z_m^n}
\left\|f\left(x+\frac{m}{2}e_j\right)-f\left(x\right)\right\|_X^pd\mu(x)
 \nonumber\\&
& \qquad-\int_{\Z_m^n}
\left\|\op_j^{(k)}f\left(x+\frac{m}{2}e_j\right)
-f\left(x+\frac{m}{2}e_j\right)\right\|_X^pd\mu(x)\nonumber\\
&&\qquad-\int_{\Z_m^n}
\left\|\op_j^{(k)}f\left(x\right)-f\left(x\right)\right\|_X^pd\mu(x)\nonumber\\
&&\quad =\frac{1}{3^{p-1}}\int_{\Z_m^n}
\left\|f\left(x+\frac{m}{2}e_j\right)-f\left(x\right)\right\|_X^pd\mu(x)\nonumber\\
&&\qquad-2\int_{\Z_m^n}
\left\|\op_j^{(k)}f\left(x\right)-f\left(x\right)\right\|_X^pd\mu(x)\nonumber\\&&\quad \ge
\frac{1}{3^{p-1}}\int_{\Z_m^n}
\left\|f\left(x+\frac{m}{2}e_j\right)-f\left(x\right)\right\|_X^pd\mu(x)\nonumber\\
&&\qquad -
2^{p+1}k^p
\E_\e\int_{\Z_m^n}\|f(x+\e)-f(x)\|_X^pd\mu(x) \nonumber\\& &\qquad - 2^{p}\int_{\Z_m^n}
\|f(x+e_j)-f(x)\|_X^pd\mu(x) ,\nonumber
\end{eqnarray} where we used Lemma~\ref{lem:approx}.

Combining~\eqref{eq:each j} with~\eqref{eq:before lemma},
we see that
\begin{align}\label{eq:almost done0}&\\[-2pt]
& \sum_{j=1}^n\int_{\Z_m^n}
\left\|f\left(x+\frac{m}{2}e_j\right)-f\left(x\right)\right\|_X^pd\mu(x)\nonumber\\
&\quad \le
\frac{(3Cm)^p n^{1-\frac{p}{q}}}{3\cdot 4^p}\E_\e
\int_{\Z_m^n}\Le\|\sum_{j=1}^n \e_j\left[\op_j^{(k)}
f(x+e_j)-\op_j^{(k)}
f(x-e_j)\right]\Ri\|_X^pd\mu(x)\nonumber\\ \nonumber  &\qquad +
6^pk^pn\E_\e\int_{\Z_m^n}\|f(x+\e)-f(x)\|_X^pd\mu(x) \\ \nonumber
&\qquad
+6^p\sum_{j=1}^n\int_{\Z_m^n}
\|f(x+e_j)-f(x)\|_X^pd\mu(x)\nonumber \displaybreak[2]\\
 &\quad \le \Le(\frac{(18Cm)^p n^{1-\frac{p}{q}}}{ 4^p}
+6^pk^pn\Ri)\E_\e\int_{\Z_m^n}\|f(x+\e)-f(x)\|_X^pd\mu(x)\nonumber\\
 &\qquad+ \Le(\frac{(3Cm)^p
n^{1-\frac{p}{q}}}{
4^p}\cdot\frac{24^pn^{2p-1}}{k^{p}}+6^p\Ri)\sum_{j=1}^n\int_{\Z_m^n}
\|f(x+e_j)-f(x)\|_X^pd\mu(x)
\nonumber\\
&\label{eq:almost done}\\[-8pt]
&\quad \le (18Cm)^p 
n^{1-\frac{p}{q}}
\left(\E_\e\int_{\Z_m^n}\|f(x+\e)-f(x)\|_X^pd\mu(x)d\sigma(\e)\right.
\nonumber\\
&\qquad \left.+\frac{1}{n}\sum_{j=1}^n\int_{\Z_m^n}
\|f(x+e_j)-f(x)\|_X^pd\mu(x)\right),\nonumber
\end{align}
where in~\eqref{eq:almost done0} we used~\eqref{eq:E},
and~\eqref{eq:almost done} holds true when we choose $4n^2\le k\le
\frac{3m}{4n^{1/q}}$ (which is possible if we assume that $m\ge
6n^{2+1/q}$). By Lemma~\ref{lem:with zeros}, this completes the
proof of Theorem~\ref{thm:weak cotype}.
\end{proof}

\subsection{Proof of Lemma~{\rm \ref{lem:pass to average on e}}}\label{section:the lemma}
 Fix $\e\in \{-1,1\}^n$, and $x\in \Z_m^n$. Consider
the following two sums:
\begin{eqnarray}\label{eq:A}
A_f(x,\e)&= &\sum_{j=1}^n
\e_j\left[\op_j^{(k)}f(x+e_j)-\op_j^{(k)}f(x-e_j)\right]\\
&=&\frac{1}{k(k+1)^{n-1}}\sum_{y\in
\Z_m^n} a_y(x,\e)f(y),\nonumber
\end{eqnarray}
and
\begin{eqnarray}\label{eq:B}
B_f(x,\e)&= &\frac{1}{k(k+1)^{n-1}} \sum_{z- x \in (-k,k)^n\cap
(2\Z)^n}[f(z+\e)-f(z-\e)]\\
&=&\frac{1}{k(k+1)^{n-1}}\sum_{y\in \Z_m^n}
b_y(x,\e)f(y),\nonumber
\end{eqnarray}
where $a_y(x,\e),b_y(x,\e)\in \Z$ are appropriately chosen
coefficients, which are independent of $f$.

For $x\in \Z_m^n$ define $S(x)\subset \Z_m^n$,
\begin{multline*}
S(x)=\Big\{y\in x+ (2\Z+1)^n:\ d_{\Z_m^n}(y,x)=k,\\ \text{and}\
|\{j:\ |y_j-x_j|\equiv k \mod m\}|\ge2\Big\}.
\end{multline*}

\begin{claim}\label{claim:cases} For $x\in \Z_m^n$ and $y\notin
S(x)${\rm ,} $a_y(x,\e)=b_y(x,\e)$.
\end{claim}
\begin{proof} If there exists a coordinate $j\in \{1,\ldots,n\}$
such that $x_j-y_j$ is even, then it follows from our definitions
that $a_y(x,\e)=b_y(x,\e)=0$. Similarly, if $d_{\Z_m^n}(x,y)>k$
then  $a_y(x,\e)=b_y(x,\e)=0$ (because $k$ is odd). Assume that
$x-y\in (2\Z+1)^n$. If $d_{\Z_m^n}(y,x)<k$ then for each $j$ the
term $f(y)$ cancels in $\op_j^{(k)}f(x+e_j)-\op_j^{(k)}(x-e_j)$,
implying that $a_y(x,\e)=0$. Similarly, in the sum defining
$B_f(x,\e)$ the term $f(y)$ appears twice, with opposite signs, so
that $b_y(x,\e)=0$.

It remains to deal with the case $|\{j:\ |y_j-x_j|\equiv k \mod
m\}|=1$. We may assume without loss of generality that
$$|y_1-x_1|\equiv k\mod m
\quad \text{ and  for $j\ge 2$},\quad  y_j-x_j\in (-k,k)\mod m.
$$
If $y_1-x_1\equiv k\mod m$ then $a_y(x,\e)=\e_1$, since in the
terms corresponding to $j\ge 2$ in the definition of $A_f(x,\e)$
the summand $f(y)$ cancels out. We also claim that in this case
$b_y(x,\e)=\e_1$. Indeed, if $\e_1=1$ then $f(y)$ appears in the
sum defining $B_f(x,\e)$ only in the term corresponding to
$z=y-\e$, while if $\e_1=-1$ then $f(y)$ appears in this sum only
in the term corresponding to $z=y+\e$, in which case its
coefficient is $-1$. In the case $y_1-x_1\equiv -k\mod m$ the same
reasoning shows that $a_y(x,\e)=b_y(x,\e)=-\e_1$.
\end{proof}

By Claim~\ref{claim:cases} we have
\begin{eqnarray}\label{eq:for later}
A_f(x,\e)-B_f(x,\e)=\frac{1}{k(k+1)^{n-1}}\sum_{y\in
S(x)}[a_y(x,\e)-b_y(x,\e)]f(y).
\end{eqnarray}

Thus,
\begin{eqnarray*}
\int_{\Z_m^n}\left\|A_f(x,\e)\right\|_X^pd\mu(x)&\le&
3^{p-1}\int_{\Z_m^n}\left\|B_f(x,\e)\right\|_X^pd\mu(x)\\
&&+3^{p-1}\int_{\Z_m^n}\Le\|\frac{1}{k(k+1)^{n-1}}\sum_{y\in
S(x)}a_y(x,\e)f(y)\Ri\|_X^pd\mu(x)
\\& &+3^{p-1}\int_{\Z_m^n}\Le\|\frac{1}{k(k+1)^{n-1}}\sum_{y\in
S(x)}b_y(x,\e)f(y)\Ri\|_X^pd\mu(x).
\end{eqnarray*}
Thus Lemma~\ref{lem:pass to average on e} will be proved once we
establish the following inequalities
\begin{equation}\label{eq:goal0}
\int_{\Z_m^n}\left\|B_f(x,\e)\right\|_X^pd\mu(x)\le
\int_{\Z_m^n}\|f(x+\e)-f(x-\e)\|_X^pd\mu(x),
\end{equation}
\begin{multline}\label{eq:goal}
\int_{\Z_m^n}\Le\|\frac{1}{k(k+1)^{n-1}}\sum_{y\in
S(x)}a_y(x,\e)f(y)\Ri\|_X^pd\mu(x)\\
\le
\frac{8^{p}n^{2p-1}}{k^p}\sum_{j=1}^n\int_{\Z_m^n}\left\|f(x+e_j)-f(x)\right\|_X^p,
\end{multline}
and
\begin{multline}\label{eq:goal2}
\int_{\Z_m^n}\Le\|\frac{1}{k(k+1)^{n-1}}\sum_{y\in
S(x)}b_y(x,\e)f(y)\Ri\|_X^pd\mu(x)\\
\le
\frac{8^{p}n^{2p-1}}{k^p}\sum_{j=1}^n\int_{\Z_m^n}\left\|f(x+e_j)-f(x)\right\|_X^p.
\end{multline}
Inequality~\eqref{eq:goal0} follows directly from the definition
of $B_f(x,\e)$, by convexity. Thus, we pass to the proof
of~\eqref{eq:goal} and~\eqref{eq:goal2}.

 For $j=1,2,\ldots,n$ define for $y\in
S(x)$,
$$
\tau_j^{x}(y)=\left\{\begin{array}{ll} y-2ke_j & y_j-x_j\equiv k\mod m,\\
y & \text{otherwise},
\end{array}\right.
$$
and set $\tau_j^{x}(y)=y$ when $y\notin S(x)$. Observe that the
following identity holds true:
\begin{eqnarray}\label{eq:invariance}
\tau_j^x(y)=\tau_j^0(y-x)+x.
\end{eqnarray}

\begin{claim}\label{claim:convexity} Assume that for every
$j\in\{1,2,\ldots,n\}${\rm ,} $x,y\in \Z_m^n$ and $\e\in \{-1,1\}^n${\rm ,} we
are given a real number $\eta_j(x,y,\e)\in [-1,1]$. Then
\begin{multline*}
\int_{\Z_m^n}\Le\|\frac{1}{k(k+1)^{n-1}}\sum_{j=1}^n\sum_{y\in
\Z_m^n}\eta_j(x,y,\e)\left[f(y)-f(\tau_j^{x}(y))\right]\Ri\|_X^pd\mu(x)\\
\le
\frac{8^{p}n^{2p-1}}{2k^p}\sum_{j=1}^n\int_{\Z_m^n}\left\|f(x+e_j)-f(x)\right\|_X^pd\mu(x).
\end{multline*}
\end{claim}

\begin{proof} Denote by $N(x,\e)$ the number of nonzero summands in
$$
\sum_{j=1}^n\sum_{y\in
\Z_m^n}\eta_j(x,y,\e)\left[f(y)-f(\tau_j^{x}(y))\right].
$$
For every $\ell\ge 2$ let $S^\ell(x)$ be the set of all $y\in
S(x)$ for which the number of coordinates $j$ such that
$y_j-x_j\in \{k,-k\}\mod m$ equals $\ell$. Then $|S^\ell(x)|=
\binom{n}{\ell}2^\ell (k-1)^{n-\ell}$. Moreover, for $y\in
S^\ell(x)$ we have that $y\neq \tau_j^x(y)$ for at most $\ell$
values of $j$. Hence
\begin{eqnarray*}
N(x,\e)\le \sum_{\ell=2}^n |S^\ell(x)|\ell&=&\sum_{\ell=2}^n
\binom{n}{\ell}2^\ell (k-1)^{n-\ell}\ell\\
&=&2n\left[
(k+1)^{n-1}-(k-1)^{n-1}\right] \le \frac{4n^2}{k^2}k(k+1)^{n-1}.
\end{eqnarray*}

Now, using \eqref{eq:invariance}, we get
\begin{eqnarray}\label{eq:step1}&&\\
&& \int_{\Z_m^n}\Le\|\frac{1}{k(k+1)^{n-1}}\sum_{j=1}^n\sum_{y\in
\Z_m^n}\eta_j(x,y,\e)\left[f(y)-f(\tau_j^{x}(y))\right]\Ri\|_X^pd\mu(x)
\nonumber\\&=&\int_{\Z_m^n}\left(\frac{N(x,\e)}{k(k+1)^{n-1}}\right)^p
\Le\|\frac{1}{N(x,\e)}\sum_{j=1}^n\sum_{y\in
\Z_m^n}\eta_j(x,y,\e)\left[f(y)-f(\tau_j^{x}(y))\right]\Ri\|_X^pd\mu(x)\nonumber\\
&\le&
\int_{\Z_m^n}\frac{N(x,\e)^{p-1}}{k^{p}(k+1)^{(n-1)p}}\sum_{j=1}^n\sum_{y\in
\Z_m^n}\left\|f(y)-f(\tau_j^{x}(y))\right\|_X^p d\mu(x)\nonumber
\\&\le&\frac{4^{p-1}n^{2p-2}}{k^{2p-1}(k+1)^{n-1}}
\sum_{j=1}^n\sum_{y\in
\Z_m^n}\int_{\Z_m^n}\left\|f(y)-f(\tau_j^{x}(y))\right\|_X^p d\mu(x)\nonumber\\
&=&
\frac{4^{p-1}n^{2p-2}}{k^{2p-1}(k+1)^{n-1}}\sum_{j=1}^n\sum_{z\in
\Z_m^n}\int_{\Z_m^n}\left\|f(z+x)-f(\tau_j^0(z)+x)\right\|_X^pd\mu(x). \nonumber
\end{eqnarray}
Consider the following set:
$$
E_j=\{z\in \Z_m^n:\ \tau_j^0(z)=z-2ke_j\}.
$$
Observe that that for every $j$,
\begin{eqnarray}\label{eq:estimate Ej}
|E_j|&=&\sum_{\ell=1}^{n-1} \binom{n-1}{\ell}2^\ell
(k-1)^{n-1-\ell}\\
&\le&(k+1)^{n-1}-(k-1)^{n-1}\le
\frac{2n}{k}(k+1)^{n-1}.\nonumber
\end{eqnarray}

Using the translation invariance of the Haar measure on $\Z_m^n$
we get that
\begin{eqnarray}\label{eq:binom} 
&& \sum_{j=1}^n\sum_{z\in
\Z_m^n}\int_{\Z_m^n}\left\|f(z+x)-f(\tau_j^0(z)+x)\right\|_X^pd\mu(x)  \\
&  =&\sum_{j=1}^n\sum_{z\in
E_j}\int_{\Z_m^n}\|f(z+x)-f(z+x-2ke_j)\|_X^pd\mu(x)\nonumber\\
&=&\sum_{j=1}^n
|E_j|\int_{\Z_m^n}\|f(w)-f(w-2ke_j)\|_X^pd\mu(w)\nonumber 
\\&\le&
\frac{2n}{k}(k+1)^{n-1}\sum_{j=1}^n\int_{\Z_m^n}
\|f(w)-f(w-2ke_j)\|_X^pd\mu(w)\nonumber\\
&\le&\label{eq:step2} \frac{2n}{k}(k+1)^{n-1}\\
&&\times\sum_{j=1}^n\int_{\Z_m^n}
\left((2k)^{p-1}\sum_{t=1}^{2k}\|f(w-(t-1)e_j) 
-f(w-te_j)\|_X^p\right)d\mu(w)\nonumber\\
&\le&2^{p+1}nk^{p-1}(k+1)^{n-1}\sum_{j=1}^n\int_{\Z_m^n}\|f(z+e_j)-f(z)\|_X^pd\mu(z) ,
\nonumber
\end{eqnarray}
where in~\eqref{eq:binom} we used~\eqref{eq:estimate Ej}.
Combining \eqref{eq:step1} and \eqref{eq:step2} completes the
proof of Claim~\ref{claim:convexity}.
\end{proof}

\medskip

By Claim~\ref{claim:convexity}, inequalities \eqref{eq:goal} and
\eqref{eq:goal2}, and hence also Lemma~\ref{lem:pass to average on
e}, will be proved once we establish the following identities:
\begin{eqnarray}\label{eq:identity}
\sum_{y\in S(x)}a_y(x,\e)f(y)=\sum_{j=1}^n\sum_{y\in
\Z_m^n}\e_j\left[f(y)-f(\tau_j^{x}(y))\right].
\end{eqnarray}
and
\begin{eqnarray}\label{eq:identity2}
\sum_{y\in S(x)}b_y(x,\e)f(y)=\sum_{j=1}^n\sum_{y\in
\Z_m^n}\delta_j(x,y,\e)\left[f(y)-f(\tau_j^{x}(y))\right],
\end{eqnarray}
for some $\delta_j(x,y,\e)\in \{-1,0,1\}$.

Identity \eqref{eq:identity} follows directly from the fact that
\eqref{eq:A} implies that for every $y\in S(x)$,
$$
a_y(x,\e)=\sum_{j:\ y_j-x_j\equiv k\mod m} \e_j-\sum_{j:\
y_j-x_j\equiv-k\mod m} \e_j.
$$

It is enough to prove identity \eqref{eq:identity2} for $x=0$,
since $b_y(x,\e)=b_{y-x}(0, \e)$. To this end we note that it
follows directly from \eqref{eq:B} that for every $y\in S(0)$
$$
b_y(0,\e)=\left\{\begin{array}{ll} 1 & \exists j \  y_j\equiv
\e_jk \mod m\text{ and } \forall \ell\
y_\ell\not\equiv-\e_jk\mod m\\
-1 & \exists j \  y_j\equiv -\e_jk \mod m\text{ and } \forall
\ell\
y_\ell\not\equiv\e_jk\mod m\\
0& \text{otherwise}.\end{array}\right.
$$
For $y\in S(0)$ define
$$
y^\ominus_j=\left\{ \begin{array}{ll} -y_j& y_j\in \{k,-k\}\mod m\\
y_j& \text{otherwise}.\end{array}\right.
$$
Since $b_y(0,\e)=-b_{y^\ominus}(0,\e)$ we get that
\begin{eqnarray}\label{eq:passtotheta}
\sum_{y\in S(0)}b_y(0,\e)f(y)=\frac12\sum_{y\in
S(0)}b_y(0,\e)\left[f(y)-f(y^\ominus)\right].
\end{eqnarray}
Define for $\ell\in \{1,\ldots,n+1\}$ a vector
$y^{\ominus_\ell}\in \Z_m^n$ by
$$
y^{\ominus_\ell}_j=\left\{ \begin{array}{ll} -y_j& j<\ell \text{ and } y_j\in \{k,-k\}\mod m\\
y_j& \text{otherwise}.\end{array}\right.
$$
Then $y^{\ominus_{n+1}}=y^\ominus$, $y^{\ominus_1}=y$ and by
\eqref{eq:passtotheta}
$$
\sum_{y\in S(0)}b_y(0,\e)f(y)= \frac12\sum_{\ell=1}^n\sum_{y\in
S(0)}b_y(0,\e)\left[f(y^{\ominus_\ell})-f(y^{\ominus_{\ell+1}})\right].
$$
Since whenever $y^{\ominus_\ell}\neq y^{\ominus_{\ell+1}}$, each
of these vectors is obtained from the other by flipping the sign
of the $\ell$-th coordinate, which is in $\{k,-k\}\mod m$, this
implies the representation \eqref{eq:identity2}. The proof of
Lemma~\ref{lem:pass to average on e} is complete.\hfill \qed

\section{A nonlinear version of the Maurey-Pisier theorem}
 
In what follows we denote by $\diag(\Z_m^n)$ the graph on $\Z_m^n$
in which $x,y\in \Z_m^n$ are adjacent if for every $i\in
\{1,\ldots,n\}$,\ $x_i-y_i\in \{\pm 1 \mod m\}$.

For technical reasons that will become clear presently, given
$\ell,n\in \mathbb N$ we denote by $\B(\M;n,\ell)$ the infimum
over $\B>0$ such that for every even $m\in\mathbb N$ and for every
$f:\Z_m^n \to \M$,
\begin{eqnarray*}
\sum_{j=1}^n\int_{\Z_m^n}d_\M\left(f\left(x+\ell
e_j\right),f(x)\right)^2d\mu(x)\le \B^2 \ell^2
n\EE_{\e}\int_{\Z_m^n}d_\M(f(x+\e),f(x))^2d\mu(x).
\end{eqnarray*}

\begin{lemma}\label{lem:a priori} For every metric space
$(\M,d_\M)${\rm ,}  every $n,a\in \mathbb N${\rm ,} every even $m,r\in \mathbb
N$ with $0\le r<m${\rm ,} and every   $f:\Z_m^n\to \M${\rm ,}
\begin{multline}\label{eq:mod}
\sum_{j=1}^n\int_{\Z_m^n}d_\M\left(f\left(x+(am+r)
e_j\right),f(x)\right)^2d\mu(x)\\* \le
\min\left\{r^2,(m-r)^2\right\}\cdot
n\EE_{\e}\int_{\Z_m^n}d_\M(f(x+\e),f(x))^2d\mu(x).
\end{multline}
In particular{\rm ,} $\B(\M;n,\ell)\le 1$ for every $n\in \mathbb N$ and
every even $\ell\in \mathbb N$.
\end{lemma}

\begin{proof}  The left-hand side of~\eqref{eq:mod} depends only on
$r$, and remains unchanged if we replace $r$ by $m-r$. We may thus
assume that $a=0$ and $r\le m-r$. Fix $x\in \Z_m^n$ and $j\in
\{1,\ldots n\}$. Observe that
$$
\left\{x+\frac{1-(-1)^k}{2}\sum_{r\neq j}e_r
+ke_j\right\}_{k=0}^{r}
$$
is a path of length $r$ joining $x$ and $x+r e_j$ in the graph
$\diag(\Z_m^n)$. Thus the distance between $x$ and $x+r e_j$ in
the graph $\diag(\Z_m^n)$ equals $r$.  If
$(x=w_0,w_1,\ldots,w_{r}=x+r e_j)$ is a geodesic joining $x$ and
$x+re_j$ in $\diag(\Z_m^n)$, then by the triangle inequality
\begin{eqnarray}\label{eq:geo}
d_\M(f(x+r e_j),f(x))^2\le
r\sum_{k=1}^{r}d_\M(f(w_k),f(w_{k-1}))^2.
\end{eqnarray}
Observe that if we sum~\eqref{eq:geo} over all geodesics joining
$x$ and $x+r e_j$ in $\diag(\Z_m^n)$, and then over all $x\in
\Z_m^n$, then in the resulting sum each edge in $\diag(\Z_m^n)$
appears the same number of times. Thus, averaging this inequality
over $x\in \Z_m^n$ we get 
\begin{eqnarray*}
\int_{\Z_m^n} d_\M(f(x+r e_j),f(x))^2d\mu(x)\le
r^2\EE_{\e}[d_\M(f(x+\e),f(x))]^2.
\end{eqnarray*}
Summing over $j=1,\ldots n$ we obtain the required result.
\end{proof}

\begin{lemma}\label{lem:sub}
For every four integers $\ell,k,s,t\in \mathbb N${\rm ,}
\[
\B\left(\M;\ell k,st\right)\le \B\left(\M;\ell,s\right)\cdot
\B\left(\M;k,t\right).
\]
\end{lemma}

\begin{proof} Let $m$ be an even integer  and take a function
$f:\Z_m^{\ell k}\to \M$. Fix $x\in \Z_m^{\ell k}$ and $\e\in
\{-1,1\}^{\ell k}$. Define $g:\Z_m^{\ell}\to \M$ by
\[
g(y)=f\Le(x+\sum_{r=1}^{k}\sum_{j=1}^{\ell} \e_{j+(r-1)\ell}\cdot
y_j\cdot  e_{j+(r-1)\ell}\Ri).
\]
By the definition of $\B\left(\M;{\ell},s\right)$,  applied to $g$,
for every $\B_1> \B\left(\M;{\ell},s\right)$,
\begin{eqnarray*}\hskip-8pt
&&\hskip-12pt \sum_{a=1}^{\ell}\int_{\Z_m^{\ell}}
d_\M\Le(f\Le(x+\sum_{r=1}^{k}\sum_{j=1}^{\ell}
\e_{j+(r-1)\ell}\cdot y_j\cdot
e_{j+(r-1)\ell}+s\sum_{r=1}^{k}\e_{a+(r-1)\ell}\cdot
e_{a+(r-1)\ell}\Ri),\\&\phantom{\le}&
f\Le(x+\sum_{r=1}^{k}\sum_{j=1}^{\ell} \e_{j+(r-1)\ell}\cdot
y_j\cdot
e_{j+(r-1)\ell}\Ri)\Ri)^2d\mu_{\Z_m^{\ell}}(y)\\
&\le& \B_1^2 s^2\ell\cdot \EE_{\delta} \int_{\Z_m^{\ell}}
d_\M\Le(f\Le(x+\sum_{r=1}^{k}\sum_{j=1}^{\ell}
\e_{j+(r-1)\ell}\cdot (y_j+\delta_j)\cdot
e_{j+(r-1)\ell}\Ri),\\&\phantom{\le}&
f\Le(x+\sum_{r=1}^{k}\sum_{j=1}^{\ell} \e_{j+(r-1)\ell}\cdot
y_j\cdot e_{j+(r-1)\ell}\Ri)\Ri)^2d\mu_{\Z_m^{\ell}}(y).
\end{eqnarray*}
Averaging this inequality over $x\in \Z_m^{{\ell k}}$ and $\e\in
\{-1,1\}^{{\ell k}}$, and using the translation invariance of the
Haar measure, we get that
\begin{multline}\label{eq:ell}
\EE_{\e}\sum_{a=1}^{\ell}\int_{\Z_m^{{\ell k}}}
d_\M\Le(f\Le(x+s\sum_{r=1}^{k}\e_{a+(r-1)\ell}\cdot
e_{a+(r-1)\ell}\Ri),f(x)\Ri)^2d\mu_{\Z_m^{{\ell k}}}(x)\\
\le \B_1^2 s^2\ell \EE_{\e}
\int_{\Z_m^{\ell k}}
d_\M\left(f\left(x+\e\right),f\left(x\right)\right)^2d\mu_{\Z_m^{{\ell k}}}(x).
\end{multline}

Next we fix $x\in \Z_m^{{\ell k}}$, $u\in \{1,\ldots, \ell\}$, and
define $h_u:\Z_m^{k}\to \M$ by
$$
h_u(y)=f\Biggl(x+{s}\sum_{r=1}^{k}y_r
 \cdot e_{u+(r-1)\ell}\Biggr).
$$
By the definition of $\B\left(\M;{k},{t}\right)$, applied to
$h_u$, for every $\B_2>\B\left(\M;{k},{t}\right)$ we have 
\begin{eqnarray*}
&& \hskip-36pt
\sum_{j=1}^{k}\int_{\Z_m^{k}}
d_\M\Biggl(f\Biggl(x+{s}\sum_{r=1}^{k}y_r \cdot
e_{u+(r-1)\ell}+{st} \cdot e_{u+(j-1)\ell}\Biggr),\\
&& \hskip100pt f\Biggl(x+{s}
\sum_{r=1}^{k}y_r \cdot
e_{u+(r-1)\ell}\Biggr)\Biggr)^2d\mu_{\Z_m^{{k}}}(y)\\
&=&
\sum_{j=1}^{k}\int_{\Z_m^{k}}
d_\M\Bigl(h_u\Bigl(y+{t}e_j\Bigr),h_u(y)\Bigr)^2
d\mu_{\Z_m^{{k}}}(y) \\
&\le&\B_2^2t^2k\cdot\EE_{\e}
\int_{\Z_m^{k}}d\left(h_u\left(y+\e\right),h_u(y)\right)^2d\mu_{\Z_m^{k}}(y)\\
&=& \B_2^2t^2k \EE_{\e} \int_{\Z_m^{k}}d_\M\Biggl(f\Biggl(x+{s}
\sum_{r=1}^{k}(y_r+\e_{u+(r-1)\ell})\cdot
e_{u+(r-1)\ell}\Biggr),\\
&& \hskip100pt f\Biggl(x+{s}\sum_{r=1}^{k}y_r \cdot
e_{u+(r-1)\ell}\Biggr)\Biggr)^2d\mu_{\Z_m^{k}}(y).
\end{eqnarray*}
Summing this inequality over $u\in \{1,\ldots,\ell\}$ and
averaging over $x\in \Z_m^{{\ell k}}$, we get,
using~\eqref{eq:ell}, that
\begin{eqnarray*}
&&\hskip-5pt \sum_{a=1}^{{\ell
k}}\int_{\Z_m^{{\ell k}}}
d_\M\left(f\left(x+{st}e_a\right),f(x)\right)^2d\mu(x)\\&&\quad  \le
\B_2^2t^2k \EE_{\e}\sum_{u=1}^{\ell} \int_{\Z_m^{{\ell k}}}
d_\M\Biggl( f\Biggl(x+{s}\sum_{r=1}^{k}\e_{u+(r-1)\ell}\cdot
e_{u+(r-1)\ell}\Biggr),f\left(x\right)\Biggr)^2
d\mu(x)\\
&&\quad \le \B_2^2t^2k\cdot\B_1^2s^2\ell\E_\e \int_{\Z_m^{{\ell
k}}}d_\M\left(f\left(x+\e\right),f\left(x\right)\right)^2d\mu(x).
\end{eqnarray*}
This implies the required result.
\end{proof}

\begin{lemma}\label{lem:use the sub} Assume that there exist integers $n_0,\ell_0>1$ such
that\break $\B(\M;n_0,\ell_0) <1$. Then there exists $0<q<\infty$ such that
for every integer $n${\rm ,}
$$
m_q^{(2)}(\M;n,3n_0)\le 2\ell_0n^{\log_{n_0}\ell_0}.
$$
In particular{\rm ,} $\Gamma_q^{(2)}(\M)<\infty$.
\end{lemma}

\begin{proof}  Let $q<\infty$ satisfy
$\B(\M,n_0,\ell_0)<n_0^{-1/q}$. Iterating Lemma~\ref{lem:sub} we
get that for every integer $k$, $\B(n_0^k,\ell_0^k)\le
n_0^{-k/q}$. Denoting $n=n_0^k$ and $m=2\ell_0^k$, this implies
that for every $f:\Z_{m}^n\to \M$,
\begin{multline*}
\sum_{j=1}^n\int_{\Z_m^n}d_\M\left(f\left(x+\frac{m}{2}
e_j\right),f(x)\right)^2d\mu(x)\\
\le\frac14 m^2n^{1-\frac{2}{q}}
\EE_{\e}\int_{\Z_m^n}d_\M(f(x+\e),f(x))^2d\mu(x).
\end{multline*}
For $f:\Z_m^{n'} \to \M$, where $n'\le n$, we define $g:\Z_m^{n'}\times \Z_m^{n-n'} \to \M$ by $g(x,y)=f(x)$. Applying the above inqeuality to $g$ we obtain,
\begin{multline*} \sum_{j=1}^{n'} \int_{\Z_m^{n'}} d_{\M} \left( f\left(x+\tfrac m2 e_j \right), f(x) \right )^2 s\mu(x) \\
\le \frac 14 m^2 n^{1-\frac 2q} \EE_{\e} \int_{\Z_m^{n'}} d{\M}(f(x+\e), f(x))^2 d\mu(x) . \end{multline*}
Hence, by 
Lemma~\ref{lem:with zeros} we deduce that
$\Gamma_q^{(2)}(\M;n_0^k,2\ell_0^k)\le 3$. For general $n$, let
$k$ be the minimal integer such that $n\le n_0^k$. By
Lemma~\ref{lem:monotone} we get that $\Gamma(\M;n,2\ell_0^k)\le
3n_0^{1-2/q}\le 3n_0$. In other words,
\vskip12pt
\hfill $
\displaystyle{m_q^{(2)}(\M;n,3n_0)\le 2\ell_0^k\le 2\ell_0n^{\log_{n_0}\ell_0}.}
$ 
\end{proof}

\begin{theorem}\label{thm:reverse cotype}
Let $n>1$ be an integer{\rm ,} $m$ an even integer{\rm ,} and $s$ an integer divisible
by $4$.
Assume that $\eta\in (0,1)$ satisfies
$8^{sn}\sqrt{\eta}<\frac12${\rm ,} and that there exists a mapping
$f:\Z_m^{n}\to \M$ such that
\begin{multline}\label{eq:condition}
\sum_{j=1}^{n}\int_{\Z_m^{n}}d_\M\left(f\left(x+s
e_j\right),f(x)\right)^2d\mu(x)
\\
> (1-\eta)s^2 n
\E_{\e}\int_{\Z_m^{n}}d_\M(f(x+\e),f(x))^2d\mu(x).
\end{multline}
Then
$$
c_\M\left(\left[s/4\right]_\infty^{n}\right)\le
1+8^{sn}\sqrt{\eta}.
$$
In particular{\rm ,} if $\B(\M;n,s)=1$ then
$c_\M\left([s/4]_\infty^{n}\right)=1$.
\end{theorem}

\begin{proof} Observe first of all that~\eqref{eq:condition} and Lemma~\ref{lem:a
priori} imply that $m\ge 2s\sqrt{1-\eta}>2s-1$, so that $m\ge 2s$.
 In what follows we will use the following numerical fact: If
$a_1,\ldots,a_r\ge 0$ and $0\le b\le \frac{1}{r}\sum_{j=1}^r a_j$,
then
\begin{eqnarray}\label{eq:b}
\sum_{j=1}^r \left(a_j-b\right)^2\le
 \sum_{j=1}^r a_j^2-rb^2.
\end{eqnarray}

For $x\in \Z_m^n$ let $\G_j^+(x)$ (resp. $\G_j^-(x)$) be the set
of all geodesics joining $x$ and $x+se_j$ (resp. $x-se_j$) in the
graph $\diag(\Z_m^n)$.
As we have seen
in the proof of Lemma~\ref{lem:a priori}, since $s$ is even, these
sets are nonempty.
Notice that if $m=2s$ then $\G_j^+(x) = \G_j^-(x)$; otherwise
$\G_j^+(x) \cap \G_j^-(x) = \emptyset$.
Denote $\G_j^{\pm}(x)=\G_j^+(x) \cup
\G_j^-(x)$, and for $\pi \in \G_j^\pm(x)$,
\[ \sgn(\pi)=\begin{cases} +1 & \text{ if } \pi \in \G_j^+(x)\\
-1 & \text{ otherwise} . \end{cases} \]

Each geodesic in $\G_j^{\pm}(x)$ has
length $s$. We write each $\pi\in \G_j^{\pm}(x)$ as a sequence of
vertices $\pi=(\pi_0=x,\pi_1,\ldots,\pi_{s}=x+ \sgn(\pi) se_j)$.
Using~\eqref{eq:b} with $a_j=d_\M(f(\pi_j),f(\pi_{j-1}))$ and
$b=\frac{1}{{s}}d_\M\left(f\left(x+se_j\right),f(x)\right)$, which
satisfy the conditions of~\eqref{eq:b} due to the triangle
inequality, we get that for each $\pi\in \G_j^\pm(x)$,
\begin{multline}\label{eq:before average on paths}
\sum_{\ell=1}^{s}\left[d_\M(f(\pi_\ell),f(\pi_{\ell-1}))-
\frac{1}{{s}}d_\M\left(f\left(x +\sgn(\pi)se_j\right),f(x)\right)\right]^2\\
\le
\sum_{k=1}^{{s}}d_\M(f(\pi_\ell),f(\pi_{\ell-1}))^2-\frac{1}{{s}}
d_\M\left(f\left(x + \sgn(\pi)se_j\right),f(x)\right)^2.
\end{multline}
By symmetry $|\G_j^+(x)|=|\G_j^-(x)|$, and this value is
independent of $x\in \Z_m^{n}$ and $j\in\{1,\ldots,n\}$. Denote
$g=|\G^\pm_j(x)|$, and observe that  $g\le 2\cdot 2^{ns}$.
Averaging~\eqref{eq:before average on paths} over all $x\in
\Z_m^{n}$ and $\pi\in \G_j^\pm(x)$, and summing over $j\in
\{1,\ldots,n\}$, we get that
\begin{eqnarray}\label{eq:average geo}
&&\frac{1}{g}\sum_{j=1}^{n}\int_{\Z_m^{n}}\sum_{\pi\in \G_j^\pm (x)}
\sum_{\ell=1}^{{s}}\Bigg[d_\M(f(\pi_\ell),f(\pi_{\ell-1})) \\
&&\hskip1.5in-
\frac{1}{{s}}d_\M\left(f\left(x+ \sgn(\pi)
se_j\right),f(x)\right)\Bigg]^2d\mu(x)\nonumber\\
\nonumber &&\quad \le
sn\E_\e
\int_{\Z_m^{n}}d_\M(f(x+\e),f(x))^2\,d\mu(x)\nonumber\\*
&&\qquad -\frac{1}{s}\sum_{j=1}^{n}\int_{\Z_m^{n}}
d_\M\left(f\left(x+se_j\right),f(x)\right)^2\,d\mu(x)\nonumber\\*
&&\quad < \eta sn \E_\e\int_{\Z_m^{n}}d_\M(f(x+\e),f(x))^2\,d\mu(x).\nonumber
\end{eqnarray}
Define $\psi:\Z_{m}^{n}\to \R$ by
\begin{multline*}
\psi(x)= 2\eta sn2^{sn}
\E_\e[d_\M(f(x+\e),f(x))^2]\\ -\sum_{j=1}^{n}\sum_{\pi\in \G_j^\pm(x)}
\sum_{\ell=1}^{s}\left[d_\M(f(\pi_\ell),f(\pi_{\ell-1}))-
\frac{1}{{s}}d_\M\left(f\left(x+ \sgn(\pi)
se_j\right),f(x)\right)\right]^2.
\end{multline*}
Inequality~\eqref{eq:average geo}, together with the bound on
$g$, implies that
$$
0<\int_{\Z_m^{n}} \psi(x)d\mu(x)=\frac{1}{(2s-1)^n}\int_{\Z_m^{n}}
\sum_{\substack{y\in \Z_m^{n}\\
d_{\Z_m^{n}}(x,y)< s}}\psi(y)d\mu(x).
$$
It follows that there exists $x^0\in \Z_m^{n}$ such that
\begin{eqnarray}\label{eq:subgrid}
&&\sum_{\substack{y\in \Z_m^{n}\\
d_{\Z_m^{n}}(x^0,y)< s}}\sum_{j=1}^{n}\sum_{\pi\in
\G_j^+(x)\bigcup \G_j^-(x)} \\&&\qquad 
\sum_{\ell=1}^{{s}}\left[d_\M(f(\pi_\ell),f(\pi_{\ell-1}))-
\frac{1}{{s}}d_\M\left(f\left(y+ \sgn(\pi) se_j\right),f(y)\right)\right]^2\nonumber\\
&&\quad
<2\eta sn2^{sn}\sum_{\substack{y\in \Z_m^{n}\\
d_{\Z_m^{n}}(x^0,y)< s}}\E_\e\left[d_\M(f(y+\e),f(y))^2\right].\nonumber
\end{eqnarray}
By scaling the metric $d_\M$ we may assume without loss of
generality that
\begin{eqnarray}\label{eq:average subgrid}
\frac{1}{(2s-1)^n}\sum_{\substack{y\in \Z_m^{n}\\
d_{\Z_m^{n}}(x^0,y)<s}}\E_\e\left[d_\M(f(y+\e),f(y))^2\right]=1.
\end{eqnarray}
It follows that there exists $y^0\in \Z_m^{n}$ satisfying
$d_{\Z_m^{n}}(x^0,y^0)< s$ such that
\begin{eqnarray}\label{eq:good point}
\E_\e\left[d_\M(f(y^0+\e),f(y^0))^2\right]\ge 1.
\end{eqnarray}
By translating the argument of $f$, and multiplying
(coordinate-wise) by an appropriate sign vector in $\{-1,1\}^{n}$,
we may assume that $y^0=0$ and all the coordinates of $x^0$ are
nonnegative. Observe that this implies that every $y\in
\{0,1,\ldots,s-1\}^{n}$ satisfies $d_{\Z_m^{n}}(x^0,y)<s$.
Thus~\eqref{eq:subgrid}, and \eqref{eq:average subgrid} imply that
for every $y\in \{0,1,\ldots,s-1\}^{n}$, every $j\in
\{1,\ldots,n\}$, every $\pi\in \G_j^\pm(y)$, and every $\ell\in
\{1,\ldots,s\}$,
\begin{multline}\label{eq:paths constant}
\left|d_\M(f(\pi_\ell),f(\pi_{\ell-1}))-
\frac{1}{{s}}d_\M\left(f\left(y+ \sgn(\pi)
se_j\right),f(y)\right)\right|\\
\le \sqrt{2\eta (2s-1)^nsn2^{sn}}\le
2^{2sn}\sqrt{\eta}.
\end{multline}

\begin{claim}\label{claim:adjacent} For every $\e,\delta\in
\{-1,1\}^{n}$ and every $x\in \Z_m^{n}${\rm ,} such that $x+\e\in
\{0,1,\ldots,s-1\}^{n}${\rm ,}
$$
\left|d_\M(f(x+\e),f(x))-d_\M(f(x+\delta),f(x))\right|\le
2\sqrt{\eta}\cdot  2^{2sn}.
$$
\end{claim}

\begin{proof}  If $\e=\delta$ then there is nothing to prove, so assume that $\e_\ell=-\delta_\ell$.
Denote $S=\{j\in \{1,\ldots,n\}:\ \e_j=-\delta_j\}$ and define
$\theta,\tau\in \{-1,1\}^{n}$ by
\begin{equation*}
\theta_j=\begin{cases} -\e_\ell & j=\ell\\ \e_j & j\in S\setminus\{\ell\} \\
 1 & j\notin S \end{cases} \qquad \mathrm {and} \qquad
\tau_j=\begin{cases} -\e_\ell & j=\ell\\ \e_j & j\in S\setminus\{\ell\} \\
 -1 & j\notin S .\end{cases}
\end{equation*}
Consider the following path $\pi$ in $\diag(\Z_m^{n})$: Start at
$x+\e\in \{0,1,\ldots,s-1\}^{n}$, go in direction $-\e$ (i.e. pass
to $x$), then go in direction $\delta$ (i.e. pass to $x+\delta$),
then go in direction $\theta$ (i.e. pass to $x+\delta+\theta$),
then go in direction $\tau$ (i.e. pass to $x+\delta+\theta+\tau$),
and repeat this process $s/4$ times. It is clear from the
construction that $\pi\in \G_\ell^{\-\e_\ell}(x+\e)$. Thus,
by~\eqref{eq:paths constant} we get that
\begin{multline*}
\left|d_\M(f(x+\e),f(x))-d_\M(f(x+\delta),f(x))\right|
\\ =\left|d_\M(f(\pi_1),f(\pi_0))-d_\M(f(\pi_2),f(\pi_1))\right|\le
2\sqrt{\eta}\cdot  2^{2sn}.
\end{multline*}
\end{proof}

\begin{corollary}\label{coro:near zero} There exists a number
$A\ge 1$ such that for every $\e\in \{-1,1\}^{n}${\rm ,}
$$
\left(1-4\sqrt{\eta}\cdot  2^{2sn}\right)A\le d_\M(f(\e),f(0))\le
\left(1+4\sqrt{\eta}\cdot 2^{2sn}\right)A.
$$
\end{corollary}
\begin{proof}
 Denote $e=\sum_{j=1}^{n}e_j=(1,1,\ldots,1)$ and take
$$
A=\left(\E_\delta
\left[d_\M(f(\delta),f(0))^2\right]\right)^{1/2}.
$$
By~\eqref{eq:good point}, $A\ge 1$. By Claim~\ref{claim:adjacent}
we know that for every $\e,\delta\in \{-1,1\}^{2^s}$,
$$
d_\M(f(\e),f(0))\le d_\M(f(e),f(0))+2\sqrt{\eta}\cdot 2^{2sn}\le
d_\M(f(\delta),f(0))+4\sqrt{\eta}\cdot 2^{2sn}.
$$
Averaging over $\delta$, and using the Cauchy-Schwartz inequality,
we get that
\begin{eqnarray*}
d_\M(f(\e),f(0))&\le& \left(\E_\delta
\left[d_\M(f(\delta),f(0))^2\right]\right)^{1/2}+4\sqrt{\eta}\cdot
2^{2sn}\\
&=& A+4\sqrt{\eta}\cdot  2^{2sn} \le
\left(1+4\sqrt{\eta}\cdot 2^{2sn}\right)A.
\end{eqnarray*}

In the reverse direction we also know that
$$
A^2=\E_\delta [d_\M(f(\delta),f(0))^2]\le
\left[d_\M(f(\e),f(0))+4\sqrt{\eta}\cdot  2^{2sn}\right]^2,
$$
which implies the required result since $A\ge 1$.
\end{proof}

\begin{claim}\label{claim:path} Denote
\begin{eqnarray}\label{eq:defV}
V=\left\{x\in \Z_m^{n}:\ \forall j\ 0\le x_j\le\frac{s}{2}\
\mathrm{and}\ x_j\ \mathrm{is\ even}\right\}.
\end{eqnarray}
Then the following assertions hold true\/{\rm :}\/
\begin{enumerate} 
\item For every $x,y\in V$ there is some $z\in \{x,y\}${\rm ,} $j\in
\{1,\ldots,n\}${\rm ,} and a path $\pi\in \G_j^+(z)$ of length $s$ which
goes through $x$ and $y$. Moreover{\rm ,} we can ensure that if
$\pi=(\pi_0,\ldots,\pi_{s})$ then for all $\ell\in
\{1,\ldots,s\}${\rm ,} $\{\pi_\ell,\pi_{\ell-1}\}\cap
\{0,\ldots,s-1\}^{n}\neq \emptyset$.

\item For every $x,y\in V${\rm ,}
$d_{\diag(\Z_m^{n})}(x,y)=d_{\Z_m^{n}}(x,y)=\|x-y\|_\infty$.
\end{enumerate}
\end{claim}

\begin{proof} Let $j\in \{1,\ldots,n\}$ be such that
$|y_j-x_j|=\|x-y\|_\infty:=  t$. Without loss of generality
$y_j\ge x_j$. We will construct a path of length $s$ in
$\G_j^+(x)$ which goes through $y$. To begin with, we define
$\e^\ell,\delta^\ell\in \{-1,1\}^{n}$ inductively on $\ell$ as follows:
\begin{eqnarray*}
\e_r^\ell&=&\begin{cases} 1 & x_r+2\sum_{k=1}^{\ell-1}(\e^k_r+\delta^k_r)<y_r\\ -1 &
x_r+2\sum_{k=1}^{\ell-1}(\e^k_r+\delta^k_r)>y_r \\
 1 & x_r+2\sum_{k=1}^{\ell-1}(\e^k_r+\delta^k_r)=y_r \end{cases} \\
\noalign{\noindent and}  
\delta_r^\ell &=&\begin{cases} 1 & x_r+2\sum_{k=1}^{\ell-1}(\e^k_r+\delta^k_r)<y_r\\ -1 &
x_r+2\sum_{k=1}^{\ell-1}(\e^k_r+\delta^k_r)>y_r \\
 -1 & x_r+2\sum_{k=1}^{\ell-1}(\e^k_r+\delta^k_r)=y_r. \end{cases}
\end{eqnarray*}
If we define $a_\ell=x+\sum_{k=1}^\ell
\e^k+\sum_{k=1}^{\ell-1}\delta^k$ and $b_\ell=x+\sum_{k=1}^\ell
\e^k+\sum_{k=1}^{\ell}\delta^k$ then the sequence
$$(x,a_1,b_1,a_2,b_2,\ldots, a_{t/2-1},b_{t/2}=y)$$ is a path of
length $t$ in $\diag(\Z_m^{n})$ joining $x$ and $y$. This proves
the second assertion above. We extend this path to a path of
length $s$ (in $\diag(\Z_m^{n})$) from $x$ to $x+se_j$ as follows.
Observe that for every $1\le \ell\le t/2$,
$\e^\ell_j=\delta^\ell_j=1$. Thus
$-\e^\ell+2e_j,-\delta^\ell+2e_j\in \{-1,1\}^n$. If we define
$c_\ell=y+\sum_{k=1}^\ell
(-\e^k+2e_j)+\sum_{k=1}^{\ell-1}(-\delta^k+2e_j)$ and
$d_\ell=y+\sum_{k=1}^\ell
(-\e^k+2e_j)+\sum_{k=1}^{\ell}(-\delta^k+2e_j)$, then
$d_{t/2}=x+2te_j$. Observe that by the definition of $V$, $2t\le
s$, and $s-2t$ is even. Thus we can continue the path from
$x+2te_j$ to $x+se_j$ by alternatively using the directions
$e_j+\sum_{\ell\neq j} e_\ell$ and $e_j-\sum_{\ell\neq j} e_\ell$.
\end{proof}

\begin{corollary}\label{coro:pass to all x}
 Assume that $x\in V$. Then
for $A$ as in Corollary~{\rm \ref{coro:near zero},} we have for all
$\e\in \{-1,1\}^{n}${\rm ,}
$$
\left(1-10\sqrt{\eta}\cdot 2^{2sn}\right)A\le
d_\M(f(x+\e),f(x))\le \left(1+10\sqrt{\eta}\cdot 2^{2sn}\right)A.
$$
\end{corollary}

\begin{proof}
By Claim~\ref{claim:path} (and its proof), there exist  $j\in
\{1,\ldots,n\}$ and $\pi\in \G_j^+(0)$ such that
$\pi_1=e=(1,\ldots,1)$ and for some $k\in \{1,\ldots,s\}$,
$\pi_k=x$. Now, by~\eqref{eq:paths constant} we have that
$$
\left|d_\M\left(f(e),f(0)\right)-d_\M\left(f(\pi_{k-1}),f(x)\right)\right|\le
2\sqrt{\eta}\cdot 2^{2sn}.
$$
Observe that since $x\in V$, $x+e\in \{0,\ldots,s-1\}^{n}$. Thus
by Claim~\ref{claim:adjacent}
\begin{eqnarray*}
&&\hskip-36pt \left|d_\M\left(f(x+\e),f(x)\right)-d_\M\left(f(e),f(0)\right)\right|
\\
&&\quad \le
\left|d_\M\left(f(e),f(0)\right)-d_\M\left(f(\pi_{k-1}),f(x)\right)\right|\\& &\qquad +
\left|d_\M\left(f(\pi_{k-1}),f(x)\right)-d_\M\left(f(x+e),f(x)\right)\right|\\
&&\qquad +
\left|d_\M\left(f(x+\e),f(x)\right)-d_\M\left(f(x+e),f(x)\right)\right|
\\&&\quad \le 6\sqrt{\eta}\cdot 2^{2sn},
\end{eqnarray*}
so that the required inequalities follow from
Corollary~\ref{coro:near zero}.
\end{proof}

\begin{corollary}\label{coro:in V}
For every  distinct $x,y\in V${\rm ,}
$$
\left(1-12\sqrt{\eta}\cdot 2^{2sn}\right)A\le
\frac{d_\M(f(x),f(y))}{\|x-y\|_\infty}\le\left(1+12\sqrt{\eta}\cdot
2^{2sn}\right)A,
$$
where $A$ is as in Corollary~{\rm \ref{coro:near zero}.}
\end{corollary}
\begin{proof}
Denote $t=\|x-y\|_\infty$; we may assume that there exists
$j\in \{1,\ldots,n\}$ such that $y_j-x_j=t$. By
Claim~\ref{claim:path} there is a path $\pi\in \G_j^+(x)$ of
length $s$ such that $\pi_{t}=y$. By~\eqref{eq:paths constant} and
Corollary~\ref{coro:pass to all x} we have for every $\ell\in
\{1,\ldots,s\}$
\begin{eqnarray*}
\left|d_\M(f(\pi_\ell),f(\pi_{\ell-1}))-\frac{1}{s}d_\M\left(f\left(x+se_j\right),f(x)\right)\right|\le
\sqrt{\eta}\cdot 2^{2sn},
\end{eqnarray*}
and
$$
\left(1-10\sqrt{\eta}\cdot 2^{2sn}\right)A\le
d_\M(f(\pi_0),f((\pi_1)) \le \left(1+10\sqrt{\eta}\cdot
2^{2sn}\right)A.
$$
Thus, for all $\ell\in \{1,\ldots,s\}$,
$$
\left(1-12\sqrt{\eta}\cdot 2^{2sn}\right)A\le
d_\M(f(\pi_\ell),f(\pi_{\ell-1}))\le \left(1+12\sqrt{\eta}\cdot
2^{2sn}\right)A.
$$
Thus
\begin{eqnarray*}
d_\M(f(x),f(y))&\le& \sum_{\ell=1}^{t}
d_\M(f(\pi_\ell),f(\pi_{\ell-1}))\le t\cdot
\left(1+12\sqrt{\eta}\cdot 2^{2sn}\right)A\\
&=&\|x-y\|_\infty\cdot
\left(1+12\sqrt{\eta}\cdot 2^{2sn}\right)A.
\end{eqnarray*}
On the other hand
\begin{eqnarray*}
d_\M(f(x),f(y))&\ge&
d_\M(f(x+se_j),f(x))-d_\M(f(x+se_j),f(y))\\
&\ge& sd_\M(f(x),f(\pi_1))-s\sqrt{\eta}\cdot 2^{2sn}-\sum_{\ell=t+1}^{s}d_\M(f(\pi_\ell),f(\pi_{\ell-1}))\\
&\ge& s\left(1-10\sqrt{\eta}\cdot
2^{2sn}\right)A-s\sqrt{\eta}\cdot
2^{2sn}\\
&&-\left(s-t\right)\left(1-12\sqrt{\eta}\cdot
2^{2sn}\right)A\\
&\ge& \|x-y\|_\infty\cdot\left(1-12\sqrt{\eta}\cdot
2^{2sn}\right)A.
\end{eqnarray*}
\vglue-20pt
\end{proof}

\medskip

This concludes the proof of Theorem~\ref{thm:reverse cotype},
since the mapping $x\mapsto x/2$ is a distortion $1$ bijection
between $(V,d_{\Z_m^n})$ and $[s/4]_\infty^n$.
\end{proof}

We are now in position to prove Theorem~\ref{thm:MPcotype}.

\begin{proof}[Proof of Theorem~{\rm \ref{thm:MPcotype}}]
We   assume  that $\Gamma_q^{(2)}(\M)=\infty$ for all
$q<\infty$. By Lemma~\ref{lem:use the sub} it follows that for
every two integers $n,s>1$, $\B(\M;n,s)=1$. Now the required
result follows from Theorem~\ref{thm:reverse cotype}.
\end{proof}

\begin{lemma}\label{lem:cotype implies no grid} Let $\M$ be a
 metric space and $K>0$. Fix $q<\infty$ and assume that
 $m:=  m_q^{(2)}(\M;n,K)<\infty$. Then
$$
c_\M\left(\Z_m^n\right)\ge \frac{n^{1/q}}{2K}.
$$
\end{lemma}

\begin{proof}   Fix  a bijection $f:\Z_m^n\to \M$. Then
\begin{eqnarray*}
\frac{nm^2}{4\|f^{-1}\|_{\Lip}^2}&\le&
\sum_{j=1}^n\int_{\Z_m^n}d_\M\left(f\left(x+\frac{m}{2}e_j\right),f(x)\right)^2d\mu(x)\\&\le&
K^2 m^2
n^{1-\frac{2}{q}}\int_{\{-1,01\}^n}\int_{\Z_m^n}d_\M(f(x+\e),f(x))^2d\mu(x)d\sigma(\e)\\&\le&
K^2 m^2 n^{1-\frac{2}{q}}\|f\|_{\Lip}^2.
\end{eqnarray*}
It follows that $ \dist(f)\ge \frac{n^{1/q}}{2K}$.
\end{proof}

\begin{corollary}\label{coro:power}
Let $\F$ be a family of metric spaces and $0<q,K,\break c<\infty$. Assume
that for all $n\in \mathbb N${\rm ,} $\Gamma_q^{(2)}(\M;n,n^c)\le K$ for
every $\M\in \F$. Then for every integer $N${\rm ,}
$$
\mathcal D_N(\F)\ge \frac{1}{2cK}\left(\frac{\log N}{\log \log
N}\right)^{1/q}.
$$
\end{corollary}

We require the following simple lemma, which shows that the
problems of embedding $[m]_\infty^n$ and $\Z_m^n$ are essentially
equivalent.

\begin{lemma}\label{lem:torus grid} The grid
$[m]_\infty^n$ embeds isometrically into $\Z_{2m}^n$. Conversely{\rm ,}
$\Z_{2m}^{n}$ embeds isometrically into $[m+1]_\infty^{2mn}$.
Moreover{\rm ,} for each $\e>0${\rm ,} $\Z_{2m}^n$ embeds with distortion
$1+6\e$ into $[m+1]^{(\lceil 1/\e\rceil+1) n}$.
\end{lemma}

\begin{proof} The first assertion follows by consideration of only
elements of $\Z_{2m}^n$ whose coordinates are at most $m-1$. Next,
the Fr\'echet embedding \[x\mapsto
(d_{\Z_{2m}}(x,0),d_{\Z_{2m}}(x,1),\ldots,d_{\Z_{2m}}(x,2m-1))\in
[m+1]_\infty^{2m},\] is an isometric embedding of $\Z_{2m}$. Thus
$\Z_{2m}^n$ embeds isometrically into $[m+1]_\infty^{2mn}$. The
final assertion is proved analogously by showing that $\Z_{2m}$
embeds with distortion $1+\e$ into $[m+1]_\infty^{\lceil
1/\e\rceil+1}$. This is done by consideration of the embedding
\begin{multline*}
x\mapsto  (d_{\Z_m}(x,0),d_{\Z_m}(x,\lfloor 2\e
m\rfloor),d_{\Z_m}(x,\lfloor 4\e m\rfloor),d_{\Z_m}(x,\lfloor 6\e
m\rfloor),\ldots
\\
\ldots ,d_{\Z_m}(x,\lfloor 2\lceil 1/\e\rceil \e
m\rfloor)),
\end{multline*}
 which is easily seen to have distortion at most
$1+6\e$.
\end{proof}

We are now in position to prove Theorem~\ref{thm:dicho}.

\begin{proof}[Proof of Theorem~{\rm \ref{thm:dicho}}]
We first prove the implication $1)\implies 2)$. Let $Z$ be the
disjoint union of all finite subsets of members of $\F$, i.e.
$$
Z= \bigsqcup\left\{\N:\ |\N|<\infty\ \mathrm{and}\ \exists\ \M\in
\F,\ \N\subseteq \M\right\}.
$$
For every $k>1$ we define a metric $d_k$ on $Z$ by
$$
d_k(x,y)=\left\{ \begin{array}{ll} \frac{d_\N(x,y)}{\diam(\N)}&
\exists\ \M\in \F,\ \exists\ \N\subseteq \M\ \mathrm{s.t.}\
|\N|<\infty\ \mathrm{and}\ x,y\in \N\\
k& \mathrm{otherwise}.
\end{array}\right.
$$
Clearly $d_k$ is a metric. Moreover, by construction, for every
$K,k>1$,
$$
q^{(2)}_{(Z,d_k)}(K)\ge q^{(2)}_\F(K).
$$
Assume for the sake of contradiction that for every $K,k>1$,
$q^{(2)}_{(Z,d_k)}(K)=\infty$. In other words, for every $q<\infty$, and
$k\geq 1$, $\Gamma^{(2)}_q(Z,d_k)=\infty$. By Lemma~\ref{lem:use
the sub} it follows that for every $k\geq 1$, and every two
integers $n,s>1$,
$$
\B\left((Z,d_k);n,s\right)=1.
$$
Theorem~\ref{thm:reverse cotype} implies that
$c_{(Z,d_k)}\left([m]_\infty^{n}\right)=1$.

By our assumption there exists a metric space $X$ such that
$c_\F(X):=  D>1$. Define a metric space $X'=X\times\{1,2\}$
via $d_{X'}((x,1),(y,1))=d_{X'}((x,2),(y,2))=d_X(x,y)$ and
$d_{X'}((x,1),(y,2))=2\diam(X)$. For large enough $s$ we have that
$c_{[2^{s-3}]_\infty^{2^s}}(X')<D$. Thus $c_{(Z,d_k)}(X')<D$ for
all $k$. Define
$$
k=\frac{4\diam(X)}{\min_{\substack{x,y\in X\\x\neq y}}d_X(x,y)}.
$$
Then there exists a bijection $f:X'\to (Z,d_k)$ with
$\dist(f)<\min\{2,D\}$. Denote $L=\|f\|_{\Lip}$.

We first claim that there exists $\M\in \F$, and a finite subset
$\N\subseteq \M$, such that $|f(X')\cap \N|\ge 2$. Indeed,
otherwise, by the definition of $d_k$, for all $x',y'\in X'$,
$d_k(f(x'),f(y'))=k$. Choosing distinct $x,y\in X$, we deduce that
$$
k=d_k(f(x,1),f(y,1))\le Ld_X(x,y)\le L\diam (X),
$$
and
\begin{eqnarray*}
k&=&d_k(f(x,1),f(y,2))
\geq \frac{L}{\dist(f)} \cdot d_{X'}((x,1),(y,2))\\
&>&
\frac{L}{2}\cdot 2\diam(X)=L\diam (X),
\end{eqnarray*}
which is a contradiction.

Thus, there exists $\M\in \F$ and a finite subset $\N\subseteq \M$
such that $|f(X')\cap \N|\ge 2$. We claim that this implies that
$f(X')\subseteq \N$. This will conclude the proof of
1)~$\Longrightarrow$~2), since the metric induced by $d_k$ on $\N$
is a re-scaling of $d_\N$, so that $X$ embeds with distortion
smaller than $D$ into $\N\subseteq \M\in \F$, which is a
contradiction of the definition of $D$.

Assume for the sake of a contradiction that there exists $x'\in
X'$ such that $f(x')\notin \N$. By our assumption there are
distinct $a',b'\in X'$ such that $f(a'),f(b')\in \N$. Now,
$$
1\ge d_k(f(a'),f((b'))
\geq \frac{L}{\dist(f)} \cdot d_{X'}(a',b') >
\frac{L}{2} \cdot \min_{\substack{u,v\in X\\u\neq v}}d_X(u,v),
$$
while
\begin{eqnarray*}
\frac{4\diam(X)}{\min_{\substack{u,v\in X\\u\neq
v}}d_X(u,v)}=k&=&d_k(f(x'),f((a'))\\
&\le& L d(x',a')\le 
L\diam(X')=2L\diam(X),
\end{eqnarray*}
which is a contradiction.

To prove the implication $2)\implies 3)$ observe that in the above
argument we have shown that there exists $k,q<\infty$ such that
$\Gamma^{(2)}_q(Z,d_k)<\infty$. It follows that for some integer
$n_0$, $\B((Z,d_k);n_0,n_0)<1$, since otherwise by
Theorem~\ref{thm:reverse cotype} we would get that $(Z,d_k)$
contains, uniformly in $n$, bi-Lipschitz copies of $[n]_\infty^n$.
Combining Lemma~\ref{lem:torus grid} and Lemma~\ref{lem:cotype
implies no grid} we arrive at a contradiction. By
Lemma~\ref{lem:use the sub}, the fact that
$\B((Z,d_k);n_0,n_0)<1$, combined with Corollary~\ref{coro:power},
implies that $\mathcal D_n(Z,d_k)=\Omega((\log n)^{\alpha})$ for
some $\alpha>0$. By the definition of $(Z,d_k)$, this implies the
required result.
\end{proof}

We end this section by proving Theorem~1.8:

\begin{proof}[Proof of Theorem~{\rm 1.8}] 
Denote $|X|=n$ and
$$
\Phi=\frac{\diam(X)}{\min_{x\neq y} d(x,y)}.
$$
Write $t=4\Phi/\e$ and let $s$ be an integer divisible by $4$ such
that $s\ge\max\{n,t\}$. Then $c_{[s]^{s}_\infty}(X)\le
1+\frac{\e}{4}$. Fix a metric space $Z$ and assume that
$c_Z(X)>1+\e$. It follows that $c_Z([s]_\infty^{s})\ge
1+\frac{\e}{2}$. By Theorem~\ref{thm:reverse cotype} we deduce
that
$$
\B(Z,s,4s)\le 1-\frac{\e^2}{2^{s^2}}.
$$
By Lemma~\ref{lem:use the sub} we have that $m_q^{(2)}(\M;n,3s)\le
8sn^{\log_{s}(4s)}$, where $q\le \frac{10^s}{\e^2}$. Thus by
Lemma~\ref{lem:cotype implies no grid} and Lemma~\ref{lem:torus
grid} we see that for any integer $n\ge 8s$,
$$
c_Z\left(\left[n^5\right]_\infty^n\right)\ge
\frac{n^{1/q}}{4s}=\frac{n^{\e^2/10^s}}{4s}.
$$
Choosing $N\approx (C\gamma)^{\frac{2^{4s}}{\e^2}}$, for an
appropriate universal constant $C$, yields the required result.
\end{proof}

\section{Applications to bi-Lipschitz, uniform, and coarse
embeddings}

Let $(\N,d_\N)$ and $(\M,d_\M)$ be metric spaces. For $f:\N\to \M$
and $t>0$ we define
$$
\Omega_f(t)=\sup\{d_\M(f(x),f(y));\ d_\N(x,y)\le t\},
$$
and
$$
\omega_f(t)=\inf\{d_\M(f(x),f(y));\ d_\N(x,y)\ge t\}.
$$
Clearly $\Omega_f$ and $\omega_f$ are nondecreasing, and for
every $x,y\in \N$,
$$
\omega_f\left(d_\N(x,y)\right)\le d_\M(f(x),f(y))\le
\Omega_f\left(d_\N(x,y)\right)
$$
 With these
definitions, $f$ is uniformly continuous if $\lim_{t\to
0}\Omega_f(t)=0$, and $f$ is a uniform embedding if $f$ is
injective and both $f$ and $f^{-1}$ are uniformly continuous.
Also, $f$ is a coarse embedding if $\Omega_f(t)<\infty$ for all
$t>0$ and $\lim_{t\to \infty} \omega_f(t)=\infty$.

\begin{lemma}\label{lem:coarse restriction} Let $(\M,d_\M)$ be a metric space{\rm ,} $n$ an integer{\rm ,}
 $\Gamma>0${\rm ,} and $0< p\le
q\le r$. Then for every function $f:\ell_r^n\to \M${\rm ,} and every
$s>0${\rm ,}
$$
n^{1/q}\omega_f(2s)\le \Gamma
m_q^{(p)}(\M;n,\Gamma) \cdot \Omega_f\Le(\frac{2\pi
sn^{1/r}}{m_q^{(p)}(\M;n,\Gamma)}\Ri).
$$
\end{lemma}

\begin{proof}   Denote $m=m_q^{(p)}(\M;n,\Gamma)$, and define
$g:\Z_m^n\to \M$ by
$$
g(x_1,\ldots,x_n)=f\Le(\sum_{j=1}^n se^{\frac{2\pi i
x_j}{m}}e_j\Ri).
$$
Then
\begin{multline*}
\int_{\{-1,0,1\}^n}\int_{\Z_m^n}
d_\M(g(x+\e),g(x))^pd\mu(x)d\sigma(\e)\\
 \le  \max_{ \e\in
\{-1,0,1\}^n}\Omega_f\Le(s\Le(\sum_{j=1}^n \left|e^{\frac{2\pi i
\e_j}{m}}-1\right|^r\Ri)^{1/r}\Ri)^p \le  \Omega_f\left(\frac{2\pi
sn^{1/r}}{m}\right)^p.
\end{multline*}
On the other hand,
\begin{eqnarray*}
\sum_{j=1}^n \int_{\Z_m^n}
d_\M\left(g\left(x+\frac{m}{2}e_j\right),g(x)\right)^pd\mu(x)\ge
n\omega_f(2s)^p.
\end{eqnarray*}
By the definition of $m_q^{(p)}(\M;n,\Gamma)$ it follows that
$$
n\omega_f(2s)^p\le \Gamma^pm^p
n^{1-\frac{p}{q}}\Omega_f\left(\frac{2\pi sn^{1/r}}{m}\right)^p,
$$
as required.
\end{proof}

\begin{corollary}\label{coro:no coarse} Let $\M$ be a metric space and assume that there exist constants $c,\Gamma>0$ such that for
infinitely many integers $n${\rm ,} $m_q^{(p)}(\M;n,\Gamma)\le
cn^{1/q}$. Then for every $r>q${\rm ,} $\ell_r$ does not uniformly or
coarsely embed into $\M$.
\end{corollary}

\begin{proof} To rule out the existence of a coarse embedding
choose $s=n^{\frac{1}{q}-\frac{1}{r}}$ in Lemma~\ref{lem:coarse
restriction}. Using Lemma~\ref{lem:lower m} we get that
\begin{equation} \label{eq:ref-coarse}
\omega_f\left(2n^{\frac{1}{q}-\frac{1}{r}}\right)\le c\Gamma
\Omega_f\left(2\pi\Gamma\right). \nonumber
\end{equation}
Since $q<r$, it follows that $\liminf_{t\to \infty}
\omega_f(t)<\infty$, so that  $f$ is not a coarse embedding.

To rule out the existence of a uniform embedding, assume that
$f:\ell_r\to X$ is invertible and $f^{-1}$ is uniformly
continuous. Then there exists $\delta>0$ such that for $x,y\in
\ell_r$, if $d_\M(f(x),f(y))< \delta$ then $\|x-y\|_r<2$. It
follows that $\omega_f(2)\ge \delta$. Choosing $s=1$ in
Lemma~\ref{lem:coarse restriction}, and using Lemma~\ref{lem:lower
m}, we get that
\begin{equation}\label{eq:ref-uniform}
0<\delta\le \omega_f(2)\le c\Gamma\Omega_f\left(2\pi\Gamma\cdot
n^{\frac{1}{r}-\frac{1}{q}} \right). \nonumber
\end{equation}
Since $r>q$ it follows that $\limsup_{t\to 0} \Omega_f(t)>0$, so
that $f$ is not uniformly continuous.
\end{proof}

The following corollary contains Theorem~\ref{thm:uniform},
Theorem~\ref{thm:uniformL_p} and Theorem~\ref{thm:coarse}.

\begin{corollary} Let $X$ be a $K$-convex
Banach space. Assume that $Y$ is a Banach space which coarsely or
uniformly embeds into $X$. Then $q_Y\le q_X$. In particular{\rm ,} for
$p,q> 0${\rm ,} $L_p$ embeds uniformly or coarsely into $L_q$ if and
only if $p\le q$ or $q\le p\le 2$.
\end{corollary}

\begin{proof} By the Maurey-Pisier theorem~\cite{MP76}, for every
$\e>0$ and every $n\in \mathbb N$, $Y$ contains a $(1+\e)$ distorted copy of
$\ell_{q_Y}^n$. By Theorem~\ref{thm:K}, since $X$ is $K$-convex,
for every $q>q_X$ there exists $\Gamma<\infty$ such that
$m_q(\M;n,\Gamma) =O\left(n^{1/q}\right)$. Thus, by the proof of
Corollary~ \ref{coro:no coarse}, if $Y$ embeds coarsely or
uniformly into $X$ then $q_Y\le q$, as required.

The fact that if $p\le q$ then $L_p$ embeds coarsely and uniformly
into $L_q$ follows from the fact that in this case $L_p$, equipped
with the metric $\|x-y\|_p^{p/q}$, embeds {\em isometrically} into
$L_q$ (for $p\le q\le 2$ this is proved in~\cite{BD-CK65}, \cite{WW75}.
For the remaining cases see Remark 5.10 in~\cite{MN04}). If $2\ge
p\ge q$ then $L_p$ is linearly isometric to a subspace of $L_q$
(see e.g.~\cite{Woj96}). It remains to prove that if $p>q$ and
$p>2$ then $L_p$ does not coarsely or uniformly embed into $L_q$.
We may assume that $q\ge 2$, since for $q\le 2$, $L_q$ embeds
coarsely and uniformly into $L_2$. But, now the required result
follows from the fact that $L_q$ is $K$ convex and $q_{L_q}=q$,
$q_{L_p}=p$ (see~\cite{MS86}).
\end{proof}

We now pass to the proof of Theorem~\ref{thm:infty grid}. Before
doing so we remark that Theorem~\ref{thm:infty grid} is almost
optimal in the following sense. The identity mapping embeds
$[m]_\infty^n$ into $\ell_q^n$ with distortion $n^{1/q}$. By the
Maurey-Pisier theorem~\cite{MP76}, $Y$ contains a copy of
$\ell_{q_Y}^n$ with distortion $1+\e$ for every $\e>0$. Thus $
c_Y([m]_\infty^n)\le n^{1/q_Y}$. Additionally, $[m]_\infty^n$ is
$m$-equivalent to an equilateral metric. Thus, if $Y$ is infinite
dimensional then $c_Y([m]_\infty^n)\le m$. It follows that
$$
c_Y([m]_\infty^n)\le \min\left\{n^{1/q_Y},m\right\}.
$$

\begin{proof}[Proof of Theorem~{\rm \ref{thm:infty grid}}] Assume that
$m$ is divisible by $4$ and
$$
m\ge \frac{2n^{1/q}}{C_q(Y)K(Y)}.
$$
By Theorem~\ref{thm:K}, for every $f:\Z_m^n \to Y$,
\begin{multline*}
\sum_{j=1}^n
\int_{\Z_m^n}\left\|f\left(x+\frac{m}{2}\right)-f(x)\right\|_Y^qd\mu(x)\\
\le
\left[15C_q(Y)K(Y)\right]^qm^q\int_{\{-1,0,1\}^n}\int_{\Z_m^n}\|f(x+\e)-f(x)\|_Y^qd\mu(x)d\sigma(\e).
\end{multline*}
Thus, assuming that $f$ is bi-Lipschitz we get that
$$
\frac{nm^q}{2^q\|f^{-1}\|_{\Lip}^q}\le
\left[15C_q(Y)K(Y)\right]^qm^q\cdot \|f\|_{\Lip}^q,
$$
i.e.
$$
\dist(f)\ge \frac{n^{1/q}}{30C_q(Y)K(Y)}.
$$
By Lemma~\ref{lem:torus grid} this shows that for $m\ge
\frac{2n^{1/q}}{C_q(Y)K(Y)}$, such that $m$ is divisible by $4$, $
c_Y([m]_\infty^n)=\Omega\left(n^{1/q}\right)$. If $m<
\frac{2n^{1/q}}{C_q(Y)K(Y)}$ then the required lower bound follows
from the fact that $[m]_\infty^n$ contains an isometric copy of
$[m_1]_\infty^{n_1}$, where $m_1$ is an integer divisible by $4$,
$m_1\ge \frac{2n_1^{1/q}}{C_q(Y)K(Y)}$, and $m_1 =\Theta(m)$,
$n_1=\Theta(m^q)$. Passing to integers $m$ which are not
necessarily divisible by $4$ is just as simple. 
\end{proof}

\begin{remark} Similar arguments yield bounds on $c_Y([m]_p^n)$,
which strengthen the bounds in~\cite{MN05-proc}.
\end{remark}

\begin{remark}\label{rem:L1}
Although $L_1$ is not $K$-convex, we can still show that
$$
c_1([m]_\infty^n)=\Theta\left(\min\left\{\sqrt{n},m\right\}\right).
$$
This is proved as follows. Assume that $f:\Z_m^n\to L_1$ is
bi-Lipschitz. If $m$ is divisible by $4$, and $m\ge \pi\sqrt{n}$,
then the fact that $L_1$, equipped with the metric
$\sqrt{\|x-y\|_1}$, is isometric to a subset of Hilbert
space~\cite{WW75}, \cite{DL97}, together with
Proposition~\ref{prop:hilbert}, shows that
\begin{multline*}
\sum_{j=1}^n
\int_{\Z_m^n}\left\|f\left(x+\frac{m}{2}\right)-f(x)\right\|_1d\mu(x)\\
\le
m^2\int_{\{-1,0,1\}^n}\int_{\Z_m^n}\|f(x+\e)-f(x)\|_1d\mu(x)d\sigma(\e).
\end{multline*}
Arguing as in the proof of Theorem~\ref{thm:infty grid}, we see
that for $m\approx \sqrt{n}$,
$c_1([m]_\infty^n)=\Omega\left(\sqrt{n}\right)$. This implies the
required result, as in the proof of Theorem~\ref{thm:infty grid}.
\end{remark}

\section{Discussion and open problems}\label{section:problems}

1.
  Perhaps the most important open problem related to the
nonlinear cotype inequality on Banach spaces is whether for every
Banach space $X$ with cotype $q<\infty$, for every $1\le p\le q$
there is a constant $\Gamma<\infty$ such that
$m_q^{(p)}(X;n,\Gamma)=O\left(n^{1/q}\right)$.
By Lemma~\ref{lem:lower m} this is best possible.
In Theorem~\ref{thm:K} we proved that this is indeed the case when
$X$ is $K$-convex, while our proof of Theorem~\ref{thm:cotype}
only gives $m_q^{(p)}(X;n,\Gamma)=O\left(n^{2+1/q}\right)$.

\smallbreak 2. $L_1$ is not $K$-convex, yet we do know that
$m_2^{(1)}(L_1;n,4)=O\left(\sqrt{n}\right)$. This follows directly
from Remark~\ref{rem:L1}, Lemma~\ref{lem:multip} and
Lemma~\ref{lem:monotone}. It would be interesting to prove the
same thing for $m_2(L_1;n,\Gamma)$.

\smallbreak 3. We conjecture that the $K$-convexity assumption in
Theorem~\ref{thm:uniform} and Theorem~\ref{thm:coarse} is not
necessary. Since $L_1$ embeds coarsely and uniformly into $L_2$,
these theorems do hold for $L_1$. It seems to be unknown whether
any Banach space with finite cotype embeds uniformly or coarsely
into a $K$-convex Banach space. The simplest space for which we do
not know the conclusion of these theorems is the Schatten trace
class $C_1$ (see~\cite{Woj96}. In~\cite{T-J74} it is shown that
this space has cotype $2$). The fact that $C_1$ does not embed
uniformly into Hilbert space follows from the results
of~\cite{AMM85}, together with~\cite{Pisier78}, \cite{Kalton85}. For more
details we refer to the discussion in~\cite{BL00} (a similar
argument works for coarse embeddings of $C_1$ into Hilbert space,
by use of~\cite{Ran04}). We remark that the arguments presented here
show that a positive solution of the first problem stated above
would yield a proof of Theorem~\ref{thm:uniform} and
Theorem~\ref{thm:coarse} without the $K$-convexity assumption.

\section{Acknowledgments} We are grateful to Keith Ball for several valuable discussions. We also thank Yuri Rabinovich for
pointing out the connection to Matou\v{s}ek's BD Ramsey theorem.

\newcommand{\name}[1]{\textsc{#1}} 
\newcommand{\bibline}{\_\_\_\_\_\_\_\_\_\,}

\end{document}